\documentclass[12pt,a4paper]{article}

\usepackage[reqno]{amsmath}
\usepackage{amssymb,amsthm,upref,amscd}
\usepackage[T1]{fontenc}
\usepackage{times}
\usepackage{color}

\newtheorem{theorem}{Theorem}[section]

\newtheorem{lemma}[theorem]{Lemma}
\newtheorem{proposition}[theorem]{Proposition}
\newtheorem{remark}[theorem]{Remark}

\theoremstyle{definition} \theoremstyle{remark}
\numberwithin{equation}{section}

\begin{document}
	
	\title{\textbf{A New Method to prove the Existence, Non-existence, Multiplicity, Uniqueness, and Orbital Stability/Instability of standing waves for NLS with partial confinement \\
		}}
		\author{Linjie Song$^{\mathrm{a,b}}$\thanks{The author is supported by CEMS. Email: songlinjie18@mails.ucas.edu.cn}, Hichem Hajaiej$^{\mathrm{c}}$\thanks{hhajaie@calstatela.edu.}\  \\
			\\
			{\small $^{\mathrm{a}}$Institute of Mathematics, AMSS, Chinese Academy of Science,
				Beijing 100190, China}\\
			{\small $^{\mathrm{b}}$University of Chinese Academy of Science,
				Beijing 100049, China}\\
			{\small $^{\mathrm{c}}$ Department of Mathematics, California State University at Los Angeles,
				Los Angeles, CA 90032, USA.}
		}
		\date{}
		\maketitle
		
		\begin{abstract}
			
			We give a new method to prove the existence, non-existence, multiplicity, orbital stability/instability of standing waves for NLS with partial confinement without the subcritical hypothesis, even in the reduction equation. Using this method, we give an affirmative answer for an open problem proposed by \cite[Remark 1.10]{BBJV} where the authors conjectured the existence of more than a normalized solution. We also establish uniqueness results of the ground state solutions depending on the bifurcation parameters. We explain that when the effect of partial confinement is strong, a dimension reduction appears for some parameters. We also find different bifurcation phenomena from the cases with full confinement.
			
			\bigskip
			
			

			\noindent\textbf{Data availability statement:} My manuscript has no associate data.
			
		\end{abstract}
		
		
		\medskip
		
		\section{Introduction}
		
		In \cite{BBJV}, Bellazzini and al studied the existence, stability, qualitative, and symmetry properties of the standing waves associated to the following Cauchy problem:
		\begin{equation} \label{eq1.1}
		\left\{
		\begin{array}{cc}
		i\partial_{t} \Psi + \Delta_{x} \Psi - (x_{1}^{2} + x_{2}^{2})\Psi + |\Psi|^{p-2}\Psi = 0, (t,x) \in \mathbb{R} \times \mathbb{R}^{3}, \\
		\Psi(0,x) = \Psi_0(x),
		\end{array}
		\right.
		\end{equation}
		with $10/3 < p < 6$. \eqref{eq1.1} models the propagation of the wave functions in Bose-Einstein condensate (BEC). BEC, the fifth state of matter, is one of the most remarkable discoveries in the twenty century. It describes a phenomenon that at very low temperatures all the atoms stop behaving like individual waves and merge together to form a dense, indistinguishable atomic wave. This phenomenon was predicted by S. N. Bose and A. Einstein in 1925, but it was first experimentally realized in dilute alkali gases in 1995 (2001 Nobel Prize in Physics attributed to E. A. Cornell, C. E Wieman, and W. Ketterle). BEC does not only provide a new tool to investigate the quantum properties of matter but also opens new perspectives for understanding the elusive phenomena of superconductivity and superfluidity. In the experiment, BEC is observed in presence of a confined potential trap and its macroscopic behavior strongly depends on the shape of this trap potential. The main goal of this paper is to provide a new approach to the existence of static solutions to
		\eqref{eq1.1}, and to provide information about the multiplicity and the stability/ instability of these solutions. The establishment of such properties is of enormous importance to the designers of BECs. For example, precise information about non-existence, or instability or multiplicity give valuable piece of information about the possibility of the realization of the BECs, and the validity of the model. Note that the physically relevant cubic nonlinearity $|\Psi|^{2}\Psi$ is covered by our study. Cubic NLS, often referred as Gross-Pitaevskii equation (GPE), gives a good description for BEC at temperature much smaller than the critical temperature, see e.g. \cite{GSS,DGPS}. (\ref{eq1.1}) covers the case of attractive interactions and partial confinement, including the limit case of the so-called cigar-shaped model, see \cite{BC}. 
		
		We emphasize that the partial confinement corresponds to the standard modelling for magnetic traps in BEC, see \cite{JP,PS} for a more detailed account. Our results provide valuable information to physicists and engineers. For example, the uniqueness of the ground state solution (Theorem \ref{thmA.3}) is the ideal situation for realization of a BEC. This finding shows that turning off the confinement in a direction does not prevent physicists/engineers from the creation of a stable wavefunction (BEC). This considerably reduces the costs. Additionally, discussions with the Photonics group of Caltech was the source of motivation to obtain the asymptotics in Theorem \ref{thmA.3}. Theorem \ref{thm1.4} has also an important physical significance, it tells us that under suitable conditions, the realization of BEC under partial confinement is only possible when the frequency of the electromagnetic field is low. 
		
		Before explaining our contribution and the challenges encountered, let us first introduce some notations: A standing wave of \eqref{eq1.1} is a solution of the form $\Psi(t,x) = e^{-i\lambda t}u(x)$ where $u$ solves $\Delta u - (x_{1}^{2} + x_{2}^{2})u + \lambda u + |u|^{p-2}u = 0$ in $\mathbb{R}^{3}$. Without loss of generality, we always assume that $u$ is real-valued. Thus, we aim to find solutions of the following equation:
		\begin{equation} \label{eq1.2}
			\left\{
			\begin{array}{cc}
				-\Delta u + (x_{1}^{2} + x_{2}^{2})u = \lambda u + |u|^{p-2}u \ in \ \mathbb{R}^{3}, \\
				u(x) \rightarrow 0 \ as \ |x| \rightarrow +\infty.
			\end{array}
			\right.
		\end{equation}
		Set
		$$
		H = \{u \in H^{1}(\mathbb{R}^{3}): \int_{\mathbb{R}^{3}}(x_{1}^{2} + x_{2}^{2})|u(x_{1},x_{2},x_{3})|^{2}dx < +\infty\},
		$$
		with the inner product
		$$
		\langle u,v \rangle = \int_{\mathbb{R}^{3}}(\nabla u \nabla v + (x_{1}^{2} + x_{2}^{2} + 1)uv)dx,
		$$
		and the norm
		$$
		\|u\|^{2} = \int_{\mathbb{R}^{3}}(|\nabla u|^{2} + (x_{1}^{2} + x_{2}^{2} + 1)u^{2})dx.
		$$
		We also introduce the corresponding functional
		$$
		J_{\lambda}(u) = \frac{1}{2}\int_{\mathbb{R}^{3}}(|\nabla u|^{2} + (x_{1}^{2} + x_{2}^{2} - \lambda)u^{2})dx - \frac{1}{p}\int_{\mathbb{R}^{3}}|u|^{p}dx
		$$
		and
		$$
		I(u) = \frac{1}{2}\int_{\mathbb{R}^{3}}(|\nabla u|^{2} + (x_{1}^{2} + x_{2}^{2}))dx - \frac{1}{p}\int_{\mathbb{R}^{3}}|u|^{p}dx.
		$$
		
		For NLS with a partial confinement, P. Antonelli, R. Carles and J. Drumond Silva, in \cite{ACS}, studied the scattering phenomena; Bellazzini, Boussaid, Jeanjean, and Visciglia, in \cite{BBJV}, gave the first result of the existence and stability of standing states in the $L^{2}$ supercritical case. In \cite{HJ}, Y.H. Hong and S.D. Jin further investigated the properties of these ground states. In \cite{Gou}, T.X. Gou studied the existence and orbital stability of standing waves to nonlinear Schr\"{o}dinger system with partial confinement in the $L^{2}$ subcritical case. Despite all these valuable contributions, many important open questions regarding \eqref{eq1.1} are still unsolved. Equation (\ref{eq1.1}) needs further understanding. First, as \cite[Remark 1.8]{BBJV} pointed, it is natural to ask whether there is a bifurcation phenomenon from the bottom of the spectrum of $\Lambda_0 = -\Delta + (x_1^2 + x_2^2)$. Note that $\Lambda_0$ is not an eigenvalue, it is unclear how to solve this problem via a standard bifurcation argument. One of our motivations is to provide an answer to this challenging question. We show that there is a $H^1(\mathbb{R}^3)$-bifurcation phenomenon from $\Lambda_0$. For $H$ norm, we will prove that this is not true. More precisely, we find a bifurcation phenomenon from infinity that is different form the case with a full confinement (see Theorem 1.4, Remark 1.6). Moreover, it is worth noticing that the power of the nonlinearity is mass-subcritical with respect to the set of variables bearing no potential, which is an essential requirement in \cite{BBJV,HJ}. Though this restriction on $p$ is not apparent in (\ref{eq1.1}) since the mass-critical exponent in one dimension is $6$. Our second motivation in this paper is to develop a new method, that can be extended to general dimensions and nonlinearities without the mass-subcritical condition. Finally, an open problem was proposed in \cite[Remark 1.10]{BBJV} whether a second solution of (\ref{eq1.1}) with small prescribed $L^{2}$ mass exists. Our last motivation is to address this important issue. We will provide an affirmative answer in Theorem $1.7$. We thank L. Jeanjean for informing us that the open problem of the existence of the second solution was already solved in \cite{WW2}. Our method is different, additionally we show that the second solution is orbitally unstable. Wei and Wu presented a novel method to study normalized solutions for combined nonlinearities when one is Sobolev critical. They then applied their method to equation (\ref{eq1.2}). On page 6 of \cite{WW2}, they stated that their method will be helpful in studying normalized solutions of other elliptic equations. We hope that this claim is true as it can give a second general approach after the one developed in \cite{Song3}. 
		
		In the previous contributions, there are two different approaches to study the existence of solutions of equation (\ref{eq1.2}) using variational methods:
		
		$\bullet$ For fixed and prescribed frequency $\lambda$, we study critical points of $J_{\lambda}$ but know nothing about the $L^{2}$ mass $\sqrt{\int_{\mathbb{R}^{3}}u^{2}dx}$.
		
		$\bullet$ For prescribed $L^{2}$ mass $c$, we study critical points of $I$ under the constraint $\int_{\mathbb{R}^{3}}u^{2}dx = c^{2}$. Then $\lambda$ is determined as
		a Lagrange multiplier and unknown.\\
		
		In \cite{BBJV}, Bellazzini and al noticed that when $2< p< 10/3,$ $ I_{c}>-\infty,$ for any $c\in \mathbb{R}$, and the classical approach by global minimization works. However, in the supercritical case $p> 10/3, I_c=-\infty$ for any $c\in \mathbb{R}$. New ideas and techniques were needed to overcome this challenging issue. In their transformative paper, \cite{BBJV}, the authors obtained static solutions by considering a suitable localized version of the minimization problem. More precisely, for every $\rho>0$, they introduced the following minimization problem
		
		\begin{equation}\label{BBJV}\tag{BBJV}
		J_r^\rho =\inf_{u\in S_c\cap B_\rho} I(u),
		\end{equation}
		where
		$$
		S_{c} = \{u \in H: \int_{\mathbb{R}^{3}}u^{2}dx = c^{2}\}, B_\rho = \{u \in H: \int_{\mathbb{R}^{3}}(|\nabla u|^{2} + (x_{1}^{2} + x_{2}^{2})u^{2})dx \leq \rho\}.
		$$
		
		They showed that for any $\rho>0$, there exists $r_0=r_0(\rho)>0$ such that
		$$
		\text{(BBJV1)}\qquad S_c\cap B_\rho\neq 0\quad\forall r<r_0
		$$
		$$
		\text{(BBJV2)}\qquad \emptyset\neq M_c^\rho \subset B_{\rho r}\quad\forall r<r_0
		$$
		where $M_c^\rho=\left\{u\in  S_c\cap B_\rho\ \text{s.t}\ I(u)=J_r^\rho\right\}.$ (BBJV1) guarantees that all the minimizing sequences of $J_r^\rho$ are compact up to the action of translations w.r.t $x_3$ provided that $r<r_0$. (BBJV2) is crucial to show that the minimizers of \eqref{BBJV} are critical points of the energy functional $I$ restricted to the sphere $S_c,$ where it is also crucial to make sure that these minimizers do not belong to the boundary of $B_{\rho}\cap S_c.$ If this holds true, then for any minimizer $u,$ there exists a Lagrange multiplier $\lambda$ such that equation \eqref{eq1.2} is satisfied.
		
		To the best of our knowledge, their method was the fist general approach to give a general line of attack to establish the existence of normalized solutions in the supercritical setting. It has then been used by many colleagues to study the minimization problems associated to various PDEs in the supercritical case. In this paper, we not only provide a simpler alternative proof of their main result, namely, Theorem $3$, \cite{BBJV}, but we also derive additional important information about the ground state solutions as well as their orbital stability/instability. In particular, we will provide  full and complete answers to some conjectures stated in \cite{BBJV}. As discussed in \cite{GSS}, the abstract framework developed by the two authors of this paper is applicable to many interesting situations. In this paper we will apply our general method to equation \eqref{eq1.2} and give an explanation on the effect of this type of confinement versus the classical one. It is worth mentioning that continuation argument which is made possible by the study of some limit problem on which the uniqueness and non-degeneracy is known was developed independently by several authors recently, see \cite{Song3,JZZ,WW2}. Let us point out that \cite{Ardilla-Carles} obtained results on the existence of ground state solutions of \eqref{eq1.2} by using a different method, see \cite[Lemma 3.2]{Ardilla-Carles}. However, their approach is heavily connected to this special case, and does not allow to answer the conjectures in \cite{BBJV}.

		\begin{theorem}\cite[Theorems 1-3]{BBJV}
			\label{thm1.1} Let $10/3 < p < 6$ and $c > 0$ sufficiently small.
			
			$(1)$ (Existence) There exists $u_c \in S_c$ such that
			$$
			I'|_{S_{c}}(u_c) = 0 \ and \ I(u_c) = \inf\{v \in H: v \in S_c, I'|_{S_{c}}(v) = 0\}.
			$$
			Moreover, $u_c$ satisfies (\ref{eq1.2}) for some $\lambda_c$.
			
			$(2)$ (Symmetry and monotonicity) $u_c$ is positive and for some $k \in \mathbb{R}$, $u_c(x_{1},x_{2},x_{3}-k)$ is radially symmetric and decreasing w.r.t. $(x_{1},x_{2})$ and w.r.t. $x_{3}$.
			
			$(3)$ (Stability) The set $M_c^\rho$ is stable under the flow associated with (\ref{eq1.1}) (the definition can be found in \cite{BBJV}).
			
			$(4)$ (Estimates of $\lambda$) For a universal constant $M > 0$ and $\Lambda_{0} = \inf \sigma(-\Delta + x_{1}^{2} + x_{2}^{2})$,
			\begin{equation}
			(1-Mc^{p-2})\Lambda_{0} \leq \lambda_c < \Lambda_{0}.
			\end{equation}
		\end{theorem}
		
		In the recent years, the study of normalized solutions (i.e. solutions with prescribed $L^{2}$ mass) has aroused great interest due to its numerous applications. There is considerable amount of valuable contributions. It is impossible to cite all of them, see e.g. \cite{Stu1, Stu2, Stu3} for mass subcritical case on $\mathbb{R}^{N}$, \cite{BN, BV, Jean} for the mass supercritical case on $\mathbb{R}^{N}$, \cite{Song, Song3} for non-autonomous nonlinearities, \cite{NTV, PV, Song2, Song3} for bounded domain cases, and \cite{Song2, Song4} for exterior domain cases. In particular, in \cite{Song3}, the authors developed a unified and very general method to study the existence of normalized solutions and orbital stability/instability of standing waves. The idea of using the behaviour of the $L^2$ norm with respect to $\lambda$ to deduce stability/instability result can also be found in \cite{LNo}. For a fixed $\lambda$, we study critical points of $J_{\lambda}$ and then give the information of the $L^{2}$ mass of these solutions when $\lambda$ changes. Inspired by \cite{Song3}, we make use of several tricks for the problem of partial confinement. For the readers' convenience, provide a short summary of these ideas here.
		
		\textbf{Step 1: } The existence of a ground state for fixed $\lambda < \Lambda_0$.
		
		Define the Nehari manifold with a parameter $\lambda$
		$$
		\mathcal{N}_{\lambda} = \{u \in H: \int_{\mathbb{R}^{3}}(|\nabla u|^{2} + (x_{1}^{2} + x_{2}^{2} - \lambda)u^{2})dx = \int_{\mathbb{R}^{3}}|u|^{p}dx\},
		$$
		and set
		$$
		h(\lambda) = \inf_{u \in \mathcal{N}_{\lambda}}J_{\lambda}(u).
		$$
		Step 1 consists in solving this minimization problem.
		
		The main difficulty in the establishment of the ground state is the lack of compactness, due to the translation invariance w.r.t. $x_3$. Using a concentration-compactness argument similar to \cite[Lemma 3.4]{BBJV}, for any $\lambda < \Lambda_{0}$ where $\Lambda_{0} = \inf \sigma(-\Delta + x_{1}^{2} + x_{2}^{2})$, we can obtain a minimizer $u_{\lambda}$ such that $J_{\lambda}(u_{\lambda}) = h(\lambda)$ (such solutions are called ground states or least action solutions in this paper). Furthermore, we will show the positivity, symmetry and monotonicity of $u_{\lambda}$.
		
		\textbf{Step 2: } The uniqueness, non-degeneracy and asymptotic behavior of $u_{\lambda}$ when $\lambda \rightarrow -\infty$.
		
		Let $\mu = 1/\lambda^{2} \rightarrow 0^+$ and
		\begin{equation} \label{eq1.4}
		v_{\mu} = |\lambda|^{-\frac{1}{p-2}}u_{\lambda}(\frac{x}{\sqrt{|\lambda}|}).
		\end{equation}
		Then $v_{\mu}$ is the ground state of the following equation
		\begin{equation} \label{eq1.5}
		\left\{
		\begin{array}{cc}
		-\Delta v + \mu(x_{1}^{2} + x_{2}^{2})v + v = |v|^{p-2}v \ in \ \mathbb{R}^{3}, \\
		v(x) \rightarrow 0 \ as \ |x| \rightarrow +\infty,
		\end{array}
		\right.
		\end{equation}
		with the corresponding functional
		$$
		\widetilde{J}_{\mu}(v) = \frac{1}{2}\int_{\mathbb{R}^{3}}(|\nabla v|^{2} + \mu(x_{1}^{2} + x_{2}^{2})v^{2} + v^{2})dx - \frac{1}{p}\int_{\mathbb{R}^{3}}|v|^{p}dx.
		$$
		We aim to prove that $v_{\mu} \rightarrow \widetilde{v}$ as $\mu \rightarrow 0^+$ where $\widetilde{v}$ is the unique (up to translations), non-degenerate, positive solution of the limit equation
		\begin{equation} \label{eq1.6}
		\left\{
		\begin{array}{cc}
		-\Delta v + v = |v|^{p-2}v \ in \ \mathbb{R}^{3}, \\
		v(x) \rightarrow 0 \ as \ |x| \rightarrow +\infty.
		\end{array}
		\right.
		\end{equation}
		This implies the uniqueness and non-degeneracy of $v_{\mu}$, which is equivalent to the uniqueness and non-degeneracy of $u_{\lambda}$.
		
		Difficulties that we will address in this step are the following:
		
		$\bullet$ The space is $H$ when $\mu > 0,$ and $H^{1}(\mathbb{R}^{3})$ when $\mu = 0$. Thus the space changes when we take the limit as $\mu \rightarrow 0^+$.
		
		$\bullet$ $\int_{\mathbb{R}^{3}}(x_{1}^{2} + x_{2}^{2})v^{2}dx$ cannot be controlled by $\|v\|_{H^{1}(\mathbb{R}^{3})}^{2}$. Therefore, it is difficult to get the boundedness of $v_{\mu}$ in $H$ before we prove the convergence of $v_{\mu}$.
		
		The first difficulty invalidates the implicit function methods if $p$ is not an integer. For this reason, we need new arguments to derive the uniqueness and non-degeneracy. To overcome the second difficulty, we show the convergence in $H^{1}(\mathbb{R}^{3})$ rather than in $H$, avoiding to show the boundedness of $\|v_{\mu}\|$. Then we will prove that $\int_{\mathbb{R}^{3}}u_{\lambda}^{2}dx \rightarrow 0,$ and standing waves associated to these ground states are orbitally unstable.
		
		\textbf{Step 3: } The uniqueness, non-degeneracy, and asymptotic behavior of $u_{\lambda}$ when $\lambda \rightarrow \Lambda_0$.
		
		Let $\tau = \Lambda_0 - \lambda \rightarrow 0^+$ and
		\begin{equation} \label{eqA.4}
		w_{\tau} = (\Lambda_0 - \lambda)^{-\frac{1}{p-2}}u_{\lambda}(x_1,x_2,\frac{x_3}{\sqrt{\Lambda_0 - \lambda}}).
		\end{equation}
		Then $w_{\tau}$ is the ground state of the following equation
		\begin{equation} \label{eqA.5}
		\left\{
		\begin{array}{cc}
		\frac{1}{\tau}(-\Delta_y + |y|^2 - \Lambda_0)w - \partial_{zz}w + w = |w|^{p-2}w \ in \ \mathbb{R}^{3}, \\
		w(x) \rightarrow 0 \ as \ |x| \rightarrow +\infty,
		\end{array}
		\right.
		\end{equation}
		where $x = (y,z) \in \mathbb{R}^2 \times \mathbb{R}$, with the corresponding functional
		$$
		\widehat{J}_{\tau}(w) = \frac{1}{2\tau}\int_{\mathbb{R}^{3}}(|\nabla_y w|^2 + |y|^2w^2 - \Lambda_0w^2)dx + \frac{1}{2}\int_{\mathbb{R}^{3}}(|\partial_z w|^{2} + w^{2})dx - \frac{1}{p}\int_{\mathbb{R}^{3}}|w|^{p}dx.
		$$
		Let
		$$
		e_1(y) = \frac{1}{\sqrt{\pi}}e^{-\frac{|y|^2}{2}}
		$$
		be the lowest eigenstate to $-\Delta_y + |y|^2$, i.e. $-\Delta_y e_1(y) + |y|^2e_1(y) = \Lambda_0e_1(y)$. When $\tau \rightarrow 0^+$, the effect of the partial confinement is so strong that  $w_{\tau} \rightarrow e_1(y)\widehat{w}(z)$ where $\widehat{w}(z)$ is the ground state of the following one-dimensional problem
		\begin{equation} \label{eqA.6}
		\left\{
		\begin{array}{cc}
		-w''(z) + w(z) = \frac{2}{p}\pi^{1-\frac{p}{2}}|w(z)|^{p-2}w(z) \ in \ \mathbb{R}, \\
		w(z) \rightarrow 0 \ as \ |z| \rightarrow +\infty,
		\end{array}
		\right.
		\end{equation}
		Then we show the uniqueness and non-degeneracy of $w_{\tau}$, which is equivalent to the uniqueness and non-degeneracy of $u_{\lambda}$, and prove that $\int_{\mathbb{R}^{3}}u_{\lambda}^{2}dx \rightarrow 0$ and standing waves associated to these ground states are orbitally stable.
		
		To summarize our previous discussion, our results read as follows:
		
		\begin{theorem}
			\label{thm1.2}  Let $2 < p < 6$. Then for $\lambda < \Lambda_{0}$, $h(\lambda)$ is achieved by a positive minimizer $u_{\lambda}$, (which is called a ground state and solves (\ref{eq1.2})). Furthermore, for some $z_\lambda \in \mathbb{R}$, the ground state $u_{\lambda}(x_{1},x_{2},x_{3}-z_\lambda)$ is radially symmetric and decreasing w.r.t. $(x_{1},x_{2})$ and w.r.t. $x_{3}$.
		\end{theorem}
		
		\begin{theorem}
			\label{thm1.3}  Let $2 < p < 6$. There exists $\Lambda_{1} < 0$ such that for $\lambda < \Lambda_{1}$, the ground state $u_{\lambda}$ of (\ref{eq1.2}), given by Theorem \ref{thm1.2}, is unique (up to translation w.r.t. $x_{3}$) in $H$ and non-degenerate in $H_s$, where
			$$
			H_s = \{u \in H: u = u(|y|,|z|), x = (y,z) \in \mathbb{R}^2 \times \mathbb{R}\}.
			$$
			Furthermore, for some $\{z_\lambda\} \subset \mathbb{R}$, as $\lambda \rightarrow -\infty$,
			$$
			|\lambda|^{-\frac{1}{p-2}}u_{\lambda}(\frac{x_1}{\sqrt{|\lambda}|},\frac{x_2}{\sqrt{|\lambda}|},\frac{x_3}{\sqrt{|\lambda}|}+z_\lambda) \rightarrow \widetilde{v} \ in \ H^1(\mathbb{R}^3),
			$$
			where $\widetilde{v}$ is the unique (up to translations), non-degenerate, positive solution of (\ref{eq1.6}).
		\end{theorem}
		
		\begin{theorem}
			\label{thmA.3}  Let $2 < p < 6$. There exists $\Lambda_{2} < \Lambda_{0}$ such that for $\Lambda_{2} < \lambda < \Lambda_{0}$, the ground state $u_{\lambda}$ of (\ref{eq1.2}), given by Theorem \ref{thm1.2}, is unique (up to translation w.r.t. $x_{3}$) and non-degenerate in $H$. Furthermore, for some $\{z_\lambda\} \subset \mathbb{R}$, as $\lambda \rightarrow \Lambda_{0}$,
			$$
			(\Lambda_{0} - \lambda)^{-\frac{1}{p-2}} u_{\lambda}(y,\frac{z}{\sqrt{\Lambda_{0} - \lambda}}+z_\lambda) \rightarrow e_1(y)\widehat{w}(z) \ in \ H,
			$$
			where $x = (y,z) \in \mathbb{R}^2 \times \mathbb{R}$ and $\widehat{w}$ is the unique (up to translations), non-degenerate, positive solution of (\ref{eqA.6}).
			
			In particular, if $10/3 < p < 6$,
			$$
			\|u_{\lambda}\|_{H^1(\mathbb{R}^3)} \rightarrow 0, \|u_{\lambda}\| \rightarrow +\infty;
			$$
			if $2 < p < 10/3$,
			$$
			\|u_{\lambda}\|_{H^1(\mathbb{R}^3)} \rightarrow 0, \|u_{\lambda}\| \rightarrow 0;
			$$
			if $p = 10/3$,
			$$
			\|u_{\lambda}\|_{H^1(\mathbb{R}^3)} \rightarrow 0, \|u_{\lambda}\| \rightarrow C_0,
			$$
			where $C_0 = \sqrt{\int_{\mathbb{R}^2}|y|^2e_1(y)^2dy\int_{\mathbb{R}}\widehat{w}(z)^2dz}$.
		\end{theorem}
		
		\begin{remark}
			Theorems \ref{thm1.3} and \ref{thmA.3} tell us that the effect of the partial confinement is so strong  when $\lambda \rightarrow \Lambda_0$ that a dimension reduction appears. On the other hand, this effect is so weak when $\lambda \rightarrow -\infty$ that we can almost ignore it. The latter is not trivial since the partial confinement is unbounded. Moreover, the dimension reduction of Bose-Einstein condensates is important both theoretically and experimentally, and it has been studied in various settings, see \cite{LSY,SY} for cigar-shaped and disk-shaped condensates, and \cite{ACM,ACS,AMSW,Bo1,Bo2,BC,BT,CH1,CH2,Fen,HT,HT} for further references and related results.
		\end{remark}
		
		\begin{remark}
			Theorem \ref{thmA.3} shows that there is a $H^1(\mathbb{R}^3)$-bifurcation phenomenon from $\Lambda_0$ for all $2 < p < 6$. For $H$-bifurcation, $\Lambda_0$ is a bifurcation point only when $2 < p < 10/3$. When $10/3 < p < 6$, there is a $H$-bifurcation phenomenon from infinity. We underline that the appearance of $10/3$ is a coincidence and it is not the mass-critical exponent for the space dimension. In fact, let us consider
			$$
			-\Delta u - |y|^\alpha u = \lambda u + |u|^{p-2}u, x = (y,z) \in \mathbb{R}^{N-d} \times \mathbb{R}^{d}, \alpha > 0, N \geq 2, 1 \leq d \leq N-1.
			$$
			Remark \ref{rmk1.5} will point that our method is applicable to this equation for all $2 < p < 2^\ast = 2N/(N - 2)^+$. The critical exponent of $H$-bifurcation is $2 + 4/(d + \alpha)$.
		\end{remark}
		
		\begin{theorem}
			\label{thm1.4}  Let $10/3 < p < 6$. Then for sufficiently small $c > 0$, (\ref{eq1.2}) admits at least two positive solutions whose $L^2$-norm are $c$ for different $\lambda$, one (for larger $\lambda$) is orbital stable while another (for smaller $\lambda$) is orbital unstable (the definition will be reviewed in Section 5).
		\end{theorem}

		\begin{theorem}
			\label{thmA.8}  Let $10/3 < p < 6$ and $u_\lambda$ be given by Theorem \ref{thm1.2}. Then
			$$
			\sup_{\lambda < \Lambda_0}\|u_\lambda\|_{L^2} \leq C,
			$$
			for some $C > 0$ independent on $\lambda$. In particular, this implies the non-existence of ground states with large $L^2$-norm.
		\end{theorem}
		
		\begin{remark}
			\label{rmk1.5}  It is worth noticing that our main results can be easily generalized. We can handle the cases of general dimension $N$ with a general nonlinearity (including mixed nonlinearity which has gained a lot of interest in the last year, see e.g. \cite{S1,S2,Ste,WW}) and a general potential ($V(|y|)$ is strictly increasing in $|y|$ and unbounded, $x = (y,z) \in \mathbb{R}^{N-d} \times \mathbb{R}^{d}, 0 \leq d \leq N$).
		\end{remark}
		
		\begin{remark}
			Consider the following equation
			\begin{equation}
			\left\{
			\begin{array}{cc}
			i\partial_{t} \Psi + \Delta_{x} \Psi - x_{1}^{2}\Psi + \Psi^3 = 0, (t,x) \in \mathbb{R} \times \mathbb{R}^{3}, \\
			\Psi(0,x) = \psi(x),
			\end{array}
			\right.
			\end{equation}
			which is the limit case of the so-called disk-shaped model. Our method yields to the existence and asymptotic behavior of the ground state solution when
			$$
			\lambda \rightarrow \inf \sigma(-\Delta + x_1^2) = \inf \sigma(-\partial_{x_1x_1} + x_1^2) = 1,
			$$
			and there is a two-dimensional reduction of Bose-Einstein condensates from three dimension. Note that this problem is mass-critical with respect to the set
			of variables bearing no potential. It needs more discussion to show the orbital stability/instability, which is an important, interesting and challenging issue.
		\end{remark}
		
		We organize this paper as follows. In Section 2, we study the existence and properties of $u_{\lambda}$ for a fixed $\lambda < \Lambda_{0}$ and give the proof of Theorem \ref{thm1.2}. In Sections 3 and 4, we discuss the uniqueness, non-degeneracy and asymptotic behavior of $u_{\lambda}$ when $\lambda \rightarrow -\infty$ and $\lambda \rightarrow \Lambda_0$ respectively, proving Theorems \ref{thm1.3} and \ref{thmA.3}. Finally in Section 5, we obtain the existence and multiplicity of normalized solutions for small $L^{2}$ norm, show their orbital stability or instability and complete the proof of Theorem \ref{thm1.4}. Moreover, in this section, we give the proof of Theorem \ref{thmA.8} and show the non-existence result.

\section{The existence and properties of $u_{\lambda}$ for a fixed $\lambda < \Lambda_{0}$}

\begin{lemma}
\label{lem2.1}	Let $\lambda < \Lambda_{0}$ and $\{u_n\} \subset H$ be a (PS)$_c$ sequence of $J_\lambda$ with $c \leq h(\lambda)$. Then there exist $u \in H$ and $\{z_n\} \subset \mathbb{R}$ such that $u_n(x_1,x_2,x_3 + z_n) \rightarrow u$ in $H$ up to a subsequence.
\end{lemma}

\textit{Proof.  } Standard arguments yield to the boundedness of $u_n$ in $H$. Since
$$
(\frac{1}{2} - \frac{1}{p})\int_{\mathbb{R}^{3}}|u_n|^pdx \rightarrow c,
$$
Thus $c \geq 0$. Then two cases will be treated.

Case 1: $c = 0$.
From
$$
(\frac{1}{2} - \frac{1}{p})\|u_n\|_\lambda^2 \rightarrow c = 0,
$$
we deduce that $u_n \rightarrow 0$ in $H$.

Case 2: $0 < c \leq h(\lambda)$.
Note that
$$
\liminf_{n \rightarrow +\infty}\int_{\mathbb{R}^{3}}|u_n|^pdx = \frac{2pc}{p-2} > 0.
$$
Then by \cite[Lemma 3.4]{BBJV}, for a sequence $\{z_n\} \subset \mathbb{R}$ and up to a subsequence, we have
$$
u_n(x_1,x_2,x_3 + z_n) \rightharpoonup u \neq 0 \ in \ H.
$$
On the one hand, $J'_\lambda(u) = 0$ and thus $u \in \mathcal{N}_\lambda$, yielding that $J_\lambda(u) \geq h(\lambda)$. On the other hand,
$$
J_\lambda(u) = (\frac{1}{2} - \frac{1}{p})\int_{\mathbb{R}^{3}}|u|^pdx \leq \liminf_{n \rightarrow +\infty}(\frac{1}{2} - \frac{1}{p})\int_{\mathbb{R}^{3}}|u_n|^pdx \rightarrow c \leq h(\lambda).
$$
Hence, $J_\lambda(u) = c = h(\lambda)$. Moreover, we have $\int_{\mathbb{R}^{3}}|u_n|^pdx \rightarrow \int_{\mathbb{R}^{3}}|u|^pdx$, implying that $\|u_n\|_\lambda \rightarrow \|u\|_\lambda$, where $\|u\|_\lambda^2 = \int_{\mathbb{R}^{3}}(|\nabla u|^2 + (x_1^2 + x_2^2)u^2 - \lambda u^2)dx$. Noticing the weak convergence of $u_n$ in $H$ and $\|\cdot\|_\lambda \sim \|\cdot\|$, we have $u_n \rightarrow u$ in $H$, which completes the proof.
\qed\vskip 5pt

\textbf{Proof of Theorem \ref{thm1.2}.} Let $u_n$ be a minimizing sequence for $J_\lambda$ at level $h(\lambda)$. By Ekeland variational principle (see \cite[Theorem 2.4]{Wi}), $u_n$ can be taken as (PS)$_{h(\lambda)}$ sequence. Then by Lemma \ref{lem2.1}, we know that $h(\lambda)$ can be achieved by a minimizer $u_{\lambda}$. Similar to the proof of \cite[Lemma 3.6]{Song}, we derive that $u_{\lambda}$ does not change sign. By using the strong maximum principle we get that $u_{\lambda} > 0$ (replacing $u_{\lambda}$ by $-u_{\lambda}$ if necessary). Then we obtain the radial symmetry and monotonicity properties w.r.t. the $(x_1,x_2)$ variables and $x_3$ variable respectively using moving planes techniques as in \cite{LN}. (These results can also be shown using the Steiner symmetrization and reflexion type arguments similar to \cite[Section 4]{BBJV}.)

\begin{remark}
We think that a second solution of (\ref{eq1.2}) can be obtained by using link theorem, which is sign-changing (more precisely, has two nodal domains), and its functional value belongs to $(h(\lambda),2h(\lambda))$. It seems quite challenging to obtain infinite solutions of (\ref{eq1.2}) for a fixed $\lambda$.
\end{remark}

\section{The uniqueness, non-degeneracy and asymptotic behavior of $u_{\lambda}$ when $\lambda \rightarrow -\infty$}

Let $u_\lambda$ be a positive ground state of (\ref{eq1.2}) with $\lambda \rightarrow -\infty$, $u_\lambda = u_\lambda(|y|,|z|)$, and $v_\mu$ be given by (\ref{eq1.4}) with $\mu = 1/\lambda^{2} \rightarrow 0^{+}$. Define
$$
\widetilde{\mathcal{N}}_{\mu} = \{v \in H: \int_{\mathbb{R}^{3}}(|\nabla v|^{2} + \mu(x_{1}^{2} + x_{2}^{2})v^{2} + v^{2})dx = \int_{\mathbb{R}^{3}}|v|^{p}dx\},
$$
and
$$
\widetilde{h}(\mu) = \inf_{v \in \widetilde{\mathcal{N}}_{\mu}}\widetilde{J}_{\mu}(v).
$$

\begin{lemma}
\label{lem3.1}	$J_\lambda(u_\lambda) = h(\lambda) \Leftrightarrow \widetilde{J}_{\mu}(v_\mu) = \widetilde{h}(\mu)$.
\end{lemma}

\textit{Proof.  } Let
$$
T(\lambda)u = |\lambda|^{-\frac{1}{p-2}}u(\frac{x}{\sqrt{|\lambda|}}).
$$
Direct computations imply that:
$$
u \in N_\lambda \Leftrightarrow T(\lambda)u \in \widetilde{\mathcal{N}}_{\mu},
$$
and
$$
\widetilde{J}_{\mu}(T(\lambda)u) = |\lambda|^{\frac{3}{2}-\frac{p}{p-2}}J_\lambda(u).
$$
Hence, we can affirm the conclusion.
\qed\vskip 5pt

\begin{lemma}
\label{lem3.2}	For any $(\mu, v) \in [0,+\infty) \times H \backslash \{0\}$, there exists a unique function $\widetilde{t}: [0,+\infty) \times H \backslash \{0\} \rightarrow (0, +\infty)$ such that
	\begin{equation}
		(\mu, t v) \in \mathcal{\widetilde{N}}, t > 0 \Leftrightarrow t = \widetilde{t}(\mu, v),
		\nonumber
	\end{equation}
	and $\widetilde{J}_{\mu}(\widetilde{t}(\mu, v)v) = \max_{t > 0}\widetilde{J}_{\mu}(tv)$, $\widetilde{t} \in C([0,+\infty) \times H \backslash \{0\}, (0, \infty))$, where $\mathcal{\widetilde{N}} := \cup_{\mu \geq 0}\mathcal{\widetilde{N}}_{\mu}$.
\end{lemma}

\textit{Proof.  } Note that
$$
(\mu, t v) \in \mathcal{\widetilde{N}} \Leftrightarrow \int_{\mathbb{R}^{3}}(|\nabla v|^{2} + \mu(x_{1}^{2} + x_{2}^{2})v^{2} + v^{2})dx = |t|^{p-2}\int_{\mathbb{R}^{3}}|v|^{p}dx.
$$
We define
$$
\widetilde{t}(\mu, v) = (\frac{\int_{\mathbb{R}^{3}}(|\nabla v|^{2} + \mu(x_{1}^{2} + x_{2}^{2})v^{2} + v^{2})dx}{\int_{\mathbb{R}^{3}}|v|^{p}dx})^{\frac{1}{p-2}}.
$$
Obviously, $\widetilde{t} \in C((-\infty,0] \times H \backslash \{0\}, (0, \infty))$.

Next we will show that $\widetilde{J}_{\mu}(\widetilde{t}(\mu, v)v) = \max_{t > 0}\widetilde{J}_{\mu}(tv)$. Set
$$
g(t) = \widetilde{J}_{\mu}(tv), t > 0.
$$
Then
\begin{eqnarray}
g'(t) &=& t\int_{\mathbb{R}^{3}}(|\nabla v|^{2} + \mu(x_{1}^{2} + x_{2}^{2})v^{2} + v^{2})dx - t^{p-1}\int_{\mathbb{R}^{3}}|v|^{p}dx \nonumber \\
&=& t\int_{\mathbb{R}^{3}}|v|^{p}dx(\widetilde{t}(\mu, v)^{p-2} - t^{p-2}).
\end{eqnarray}
Thus we derive that $g'(t) > 0$ if $0 < t < \widetilde{t}(\mu, v)$ and $g'(t) < 0$ if $t > \widetilde{t}(\mu, v)$. This completes the proof.
\qed\vskip 5pt

\begin{lemma}
\label{lem3.3}	$\lim_{\mu \rightarrow 0^{+}}\widetilde{h}(\mu) = \widetilde{h}(0)$.
\end{lemma}

\textit{Proof.  } Step 1: $\liminf_{\mu \rightarrow 0^{+}}\widetilde{h}(\mu) \geq \widetilde{h}(0)$.

By Lemma \ref{lem3.2}, we have
$$
\widetilde{h}(\mu) = \inf_{\|v\| = 1}\sup_{t > 0}\widetilde{J}_\mu(tv).
$$
From $\mu \geq 0$, we deduce that
$$
\sup_{t > 0}\widetilde{J}_\mu(tv) \geq \sup_{t > 0}\widetilde{J}_0(tv) \geq \inf_{\|v\| = 1}\sup_{t > 0}\widetilde{J}_0(tv) = \widetilde{h}(0), \forall v, \|v\| = 1.
$$
Hence, $\widetilde{h}(\mu) \geq \widetilde{h}(0)$, yielding that $\liminf_{\mu \rightarrow 0^{+}}\widetilde{h}(\mu) \geq \widetilde{h}(0)$.

Step 2: $\limsup_{\mu \rightarrow 0^{+}}\widetilde{h}(\mu) \leq \widetilde{h}(0)$.

Recall that $\widetilde{v}$ is the unique (up to translations), non-degenerate, positive solution of (\ref{eq1.6}). Furthermore, $\widetilde{v}$ is a ground state solution and decays exponentially at $\infty$. In particular, $\widetilde{v} \in H$. We consider $\widetilde{t}(\mu, \widetilde{v})\widetilde{v} \in \mathcal{\widetilde{N}}_{\mu}$. By Lemma \ref{lem3.2}, $\widetilde{t}(\mu, \widetilde{v})\widetilde{v} \rightarrow \widetilde{t}(0,\widetilde{v})\widetilde{v} = \widetilde{v}$ in $H$ as $\mu \rightarrow 0^{+}$, implying that:
\begin{eqnarray}
	\widetilde{h}(\mu) &\leq& \widetilde{J}_{\mu}(\widetilde{t}(\mu, \widetilde{v})\widetilde{v}) \nonumber \\
	&=& \widetilde{J}_{0}(\widetilde{t}(\mu, \widetilde{v})\widetilde{v}) + \frac{\widetilde{t}(\mu, \widetilde{v})^2\mu}{2}\int_{\mathbb{R}^{3}}(x_1^2 + x_2^2)\widetilde{v}^2dx \nonumber \\
	&=& \widetilde{J}_{0}(\widetilde{v}) + o(1) \nonumber \\
	&=& \widetilde{h}(0) + o(1).
\end{eqnarray}
Letting $\mu \rightarrow 0^{+}$, we have $\limsup_{\mu \rightarrow 0^{+}}\widetilde{h}(\mu) \leq \widetilde{h}(0)$.
\qed\vskip 5pt

\begin{lemma}
\label{lem3.4}	$v_\mu$ is bounded in $H^1(\mathbb{R}^3)$ and $\mu\int_{\mathbb{R}^{3}}(x_1^2+x_2^2)v_\mu^2dx$ is bounded as $\mu \rightarrow 0^{+}$.
\end{lemma}

\textit{Proof.  } By Lemma \ref{lem3.3}, we may assume that, as $\mu \rightarrow 0^{+}$,
\begin{equation} \label{eq3.1}
	\frac{1}{2}\int_{\mathbb{R}^{3}}(|\nabla v_\mu|^{2} + \mu(x_{1}^{2} + x_{2}^{2})v_\mu^{2} + v_\mu^{2})dx - \frac{1}{p}\int_{\mathbb{R}^{3}}|v_\mu|^{p}dx \leq \widetilde{h}(0) + 1.
\end{equation}
(\ref{eq3.1}), together with
$$
\int_{\mathbb{R}^{3}}(|\nabla v_\mu|^{2} + \mu(x_{1}^{2} + x_{2}^{2})v_\mu^{2} + v_\mu^{2})dx = \int_{\mathbb{R}^{3}}|v_\mu|^{p}dx,
$$
imply that
$$
(\frac{1}{2} - \frac{1}{p})\int_{\mathbb{R}^{3}}(|\nabla v_\mu|^{2} + \mu(x_{1}^{2} + x_{2}^{2})v_\mu^{2} + v_\mu^{2})dx \leq \widetilde{h}(0) + 1,
$$
which completes the proof.

\vskip 0.7cm

Define
$$
\widetilde{L}_\mu = -\Delta + \mu(x_1^2 + x_2^2) + 1 - (p-1)|v_\mu|^{p-2},
$$
which is a self-adjoint operator from $L^2(\mathbb{R}^3)$ to $L^2(\mathbb{R}^3)$, with form domain $H$ if $\mu > 0$ and $H^1(\mathbb{R}^3)$ if $\mu = 0$. Let $\widetilde{L}_{\mu,s}$ be the restriction of $\widetilde{L}_\mu$ on
$$
L_s^2(\mathbb{R}^3) = \{v \in L^2(\mathbb{R}^3): v = v(|y|,|z|), x = (y,z) \in \mathbb{R}^2 \times \mathbb{R}\}.
$$
We say that $v_\mu$ is non-degenerate in $H_s = H \cap L_s^2$ if $\ker \widetilde{L}_{\mu,s} = 0$.

\begin{lemma}
\label{lem3.5}  Let $v_\mu \rightarrow v_0$ in $H^1(\mathbb{R}^3)$ as $\mu \rightarrow 0^{+}$, where $v_\mu \in H_s$ solves (\ref{eq1.5}), and $v_0$ is a non-degenerate solution of (\ref{eq1.6}). Then there exists $\delta > 0$ such that $v_\mu$ is the unique and non-degenerate solution of (\ref{eq1.5}) in $\{v \in H_s: \|v - v_0\|_{H^1(\mathbb{R}^3)} < \delta\}$.
\end{lemma}

\textit{Proof.  } Step 1: Uniqueness.

We argue by contradiction and assume that there exist $v_{1,n} \neq v_{2,n}$ such that $v_{1,n} \rightarrow v_0, v_{2,n} \rightarrow v_0$ in $H^1(\mathbb{R}^3)$, $\mu_n \rightarrow 0^{+}$. Set
$$
w_n = \frac{v_{1,n} - v_{2,n}}{\|v_{1,n} - v_{2,n}\|_{H^1(\mathbb{R}^3)}}.
$$
Up to a subsequence, we may assume that $w_n \rightharpoonup w_0$ in $H^1(\mathbb{R}^3)$ and $w_n \rightarrow w_0$ in $L^p_{loc}(\mathbb{R}^3)$. Note that $w_n$ satisfies
\begin{equation} \label{eq3.3}
-\Delta w_n + \mu_n(x_1^2 + x_2^2)w_n + w_n = (p-1)|v_{2,n} + \theta_n(v_{1,n} - v_{2,n})|^{p-2}w_n, \theta_n \in [0,1].
\end{equation}
For any $\phi \in C_0^\infty(\mathbb{R}^3)$, we have
$$
\int_{\mathbb{R}^3}\mu_n(x_1^2 + x_2^2)w_n\phi dx = \int_{\mathrm{supp}\ \phi}\mu_n(x_1^2 + x_2^2)w_n\phi dx \rightarrow 0.
$$
This, together with the weak convergence of $w_n$ and the strong convergence of $v_{1,n}, v_{2,n}$ in $H^1(\mathbb{R}^3)$, we derive that
$$
\int_{\mathbb{R}^3}(\nabla w_0 \nabla \phi + w_0 \phi)dx = (p-1)\int_{\mathbb{R}^3}|v_0|^{p-2}w_0 \phi dx, \forall \phi \in C_0^\infty(\mathbb{R}^3).
$$
Since $C_0^\infty(\mathbb{R}^3)$ is dense in $H^1(\mathbb{R}^3)$, we have
$$
\int_{\mathbb{R}^3}(\nabla w_0 \nabla \phi + w_0 \phi)dx = (p-1)\int_{\mathbb{R}^3}|v_0|^{p-2}w_0 \phi dx, \forall \phi \in H^1(\mathbb{R}^3).
$$
From the non-degeneracy of $v_0$, we derive that $w_0 = 0$. Indeed, $w_0(-x) = -w_0(x)$ since $w_0 \in \mathrm{span}\{\partial_1v_0, \partial_2v_0, \partial_3v_0\}$. Noticing that $w_0(-x) = w_0(x)$ also holds, the proof is completed.

Standard arguments about (\ref{eq3.3}) yield that $w_n$ exponentially decays uniformly with respect to $n$. This, together with $w_n \rightarrow 0$ in $L^p_{loc}(\mathbb{R}^3)$, we can deduce that $w_n \rightarrow 0$ in $L^p(\mathbb{R}^3)$. Hence, by (\ref{eq3.3}),
\begin{eqnarray}
  \int_{\mathbb{R}^3}(|\nabla w_n|^2 + w_n^2)dx &\leq& \int_{\mathbb{R}^3}(|\nabla w_n|^2 + \mu_n(x_1^2 + x_2^2)w_n^2 + w_n^2)dx \nonumber \\
  &=& (p-1)\int_{\mathbb{R}^3}|v_{2,n} + \theta_n(v_{1,n} - v_{2,n})|^{p-2}w_n^2 dx \nonumber \\
  &\leq& (p-1)\|v_{2,n} + \theta_n(v_{1,n} - v_{2,n})\|_{L^p(\mathbb{R}^3)}^{p-2}\|w_n\|_{L^p(\mathbb{R}^3)}^{2} \rightarrow 0,
\end{eqnarray}
showing that $w_n \rightarrow 0$ in $H^1(\mathbb{R}^3)$, which contradicts the fact that $\|w_n\|_{H^1(\mathbb{R}^3)} = 1$.

Step 2: Non-degeneracy.

Suppose on the contrary that there exists $\varphi_n \in H_s$ satisfying
\begin{equation} \label{eq3.5}
-\Delta \varphi_n + \mu_n(x_1^2 + x_2^2)\varphi_n + \varphi_n = (p-1)|v_{\mu_n}|^{p-2}\varphi_n,
\end{equation}
where $\mu_n \rightarrow 0^+$ and $\|\varphi_n\|_{H^1(\mathbb{R}^3)} = 1$. Then similar to the proof of Step 1, we can show that $\varphi_n \rightarrow 0$ in $H^1(\mathbb{R}^3)$, in a contradiction with $\|\varphi_n\|_{H^1(\mathbb{R}^3)} = 1$.
\qed\vskip 5pt

\begin{lemma}
\label{lem3.6}  $v_\mu \rightarrow \widetilde{v}$ in $H^1(\mathbb{R}^3)$ as $\mu \rightarrow 0^{+}$.
\end{lemma}

\textit{Proof.  } Let $v_n = v_{\mu_n}$ with $\mu_n \rightarrow 0^{+}$. We will show that $v_n \rightarrow \widetilde{v}$ in $H^1(\mathbb{R}^3)$ passing to a subsequence if necessary. Lemma \ref{lem3.4} yields the boundedness of $v_n$ in $H^1(\mathbb{R}^{3})$. Then two cases will be addressed:

First, we assume that $\int_{\mathbb{R}^3}|v_n|^p dx \rightarrow 0$. Then
$$
\|v_n\|_{H^1(\mathbb{R}^3)}^2 + \mu_n\int_{\mathbb{R}^{3}}(x_1^2+x_2^2)v_n^2dx = \int_{\mathbb{R}^3}|v_n|^p dx \rightarrow 0,
$$
implying that $v_n \rightarrow 0$ in $H^1(\mathbb{R}^3)$.

Second, we assume that $\liminf_{n \rightarrow +\infty}\int_{\mathbb{R}^3}|v_n|^p dx = \alpha > 0$. Noticing that $v_n(0) = \max v_n(x)$ and that $v_n = v_n(|y|,|z|)$ is decreasing in $|y| > 0, |z| > 0$, standard methods imply that, passing to a subsequence if necessary, there exists $v_0 \in H^1(\mathbb{R}^{3})$ such that $v_n(x) \rightharpoonup v_0 \neq 0$ in $H^1(\mathbb{R}^{3})$, $v_n(x) \rightarrow v_0$ in $L^q_{loc}(\mathbb{R}^{3}), 2 \leq q < 2^\ast$, $v_n(x) \rightarrow v_0$ a.e. $x \in \mathbb{R}^{3}$.

By Lemma \ref{lem3.4}
\begin{eqnarray}
	&& |\left\langle \widetilde{J}'_0(v_n), \phi \right\rangle| \nonumber \\
	&=& \mu_n|\int_{\mathbb{R}^{3}}(x_1^2+x_2^2)v_n\phi dx| \nonumber \\
	&\leq& \sqrt{\mu_n\int_{\mathbb{R}^{3}}(x_1^2+x_2^2)v_n^2dx}\sqrt{\mu_n\int_{\mathbb{R}^{3}}(x_1^2+x_2^2)\phi^2dx} \nonumber \\
	&\rightarrow& 0, \forall \phi \in C_0^\infty(\mathbb{R}^{3}).
\end{eqnarray}
Since $C_0^\infty(\mathbb{R}^{3})$ is dense in $H^1(\mathbb{R}^3)$, we have
$$
\left\langle \widetilde{J}'_0(v_n), \phi \right\rangle \rightarrow 0, \forall \phi \in H^1(\mathbb{R}^3),
$$
implying that $\widetilde{J}'_0(v_0) = 0$. Thus $v_0 \in \widetilde{\mathcal{N}}_{0},$ and we have
\begin{eqnarray}
\widetilde{h}(0) &\leq& \widetilde{J}_0(v_0) \\ \nonumber
&=& (\frac{1}{2} - \frac{1}{p})\int_{\mathbb{R}^3}|v_0|^p dx \\ \nonumber
&\leq& \liminf_{n \rightarrow +\infty}(\frac{1}{2} - \frac{1}{p})\int_{\mathbb{R}^3}|v_n|^p dx \\ \nonumber
&\leq& \limsup_{n \rightarrow +\infty}(\frac{1}{2} - \frac{1}{p})\int_{\mathbb{R}^3}|v_n|^p dx \\ \nonumber
&=& \lim_{n \rightarrow +\infty}\widetilde{J}_{\mu_n}(v_n) \\ \nonumber
&=& \widetilde{h}(\mu_n).
\end{eqnarray}
By Lemma \ref{lem3.3}, $\widetilde{h}(\mu_n) \rightarrow \widetilde{h}(0)$. Hence, $\widetilde{J}_0(v_0) = \widetilde{h}(0)$. Moreover, we have $\int_{\mathbb{R}^3}|v_n|^p dx \rightarrow \int_{\mathbb{R}^3}|v_0|^p dx$, and then
\begin{eqnarray}
\|v_0\|_{H^1(\mathbb{R}^3)}^2 &\leq& \liminf_{n \rightarrow +\infty}\|v_n\|_{H^1(\mathbb{R}^3)}^2 \\ \nonumber
&\leq& \limsup_{n \rightarrow +\infty}\|v_n\|_{H^1(\mathbb{R}^3)}^2 \\ \nonumber
&\leq& \limsup_{n \rightarrow +\infty}(\|v_n\|_{H^1(\mathbb{R}^3)}^2 + \mu_n\int_{\mathbb{R}^{3}}(x_1^2+x_2^2)v_n^2dx) \\ \nonumber
&=& \lim_{n \rightarrow +\infty}\int_{\mathbb{R}^3}|v_n|^p dx \\ \nonumber
&=& \int_{\mathbb{R}^3}|v_0|^p dx \\ \nonumber
&=& \|v_0\|_{H^1(\mathbb{R}^3)}^2.
\end{eqnarray}
Concerning the weak convergence of $v_n$ in $H^1(\mathbb{R}^3)$, we obtain that $v_n \rightarrow v_0$ in $H^1(\mathbb{R}^3)$.

Combining the first and the second case, we prove that $v_n \rightarrow v_0$ in $H^1(\mathbb{R}^3)$ passing to a subsequence if necessary, where $v_0 = 0$ or $v_0 = \widetilde{v}$. We claim that $v_0 \neq 0$. If this is not the case, $0$ is non-degenerate ($\widetilde{J}_0''(0) = -\Delta + 1$) and thus Lemma \ref{lem3.5} yields that the solution of (\ref{eq1.6}) is unique in $\{v \in H_s: \|v\|_{H^1(\mathbb{R}^3)} < \delta\}$. Noticing that $0$ is always a solution of (\ref{eq1.6}) for any $\mu \in \mathbb{R}$, we know that $v_n = 0$. This contradicts the nontriviality of $v_n$. The claim has been proven, yielding that $v_0 = \widetilde{v}$ and the proof is completed.
\qed\vskip 5pt

\textbf{Proof of Theorem \ref{thm1.3}. } Since $\widetilde{v}$ is non-degenerate, the combination of Lemma \ref{lem3.5} and Lemma \ref{lem3.6} enables us to prove that there exists $\delta > 0$ small enough such that $v_\mu$ is the unique and non-degenerate ground state solution of (\ref{eq1.5}) in the neighborhood $\{(v,\mu) \in H_s \times [0,+\infty):\|v - \widetilde{v}\|_{H^1(\mathbb{R}^3)} + |\mu| < \delta\}$.

By Lemma \ref{lem3.1}, the uniqueness of $v_{\mu}$ is equivalent to the one of $u_{\lambda}$. Similarly, the non-degeneracy of $v_{\mu}$ is equivalent to the one of $u_{\lambda}$. Furthermore, by Theorem \ref{thm1.2}, any ground state $u_\lambda$ of (\ref{eq1.2}) satisfying $u_\lambda(0) = \max u_\lambda(x)$ belongs to $H_s$. Therefore, there exists $\Lambda_{1} < \Lambda_0$ such that when $\lambda < \Lambda_{1}$, the positive ground state of (\ref{eq1.2}) is unique in $H$ and non-degenerate in $H_s$.
\qed\vskip 5pt

\begin{remark}
Let $V = V(|y|,|z|)$ be a K-R potential, $L_\bot^2$ be the orthogonal complement of $L_s^2(\mathbb{R}^3)$ in $L^2(\mathbb{R}^3)$, and $H_{|L_\bot^2}$ be the restriction of $-\Delta + |y|^2 + V(|y|,|z|)$ to $L_\bot^2$. It is open whether $H_{|L_\bot^2}$ enjoys the Perron-Frobenius property, i.e the first eigenvalue (if it exists) of $H_{|L_\bot^2}$ is simple and the corresponding eigenfunction $\psi = \psi(|y|,|z|)$ satisfies $\psi(y,z) > 0$ for $z > 0$. If this is true, then we can get the non-degeneracy of $u_{\lambda}$ in $H$.
\end{remark}

\section{The uniqueness, non-degeneracy and asymptotic behavior of $u_{\lambda}$ when $\lambda \rightarrow \Lambda_0$}

Let $u_\lambda$ be a positive ground state of (\ref{eq1.2}) with $\lambda \rightarrow \Lambda_0$, $u_\lambda = u_\lambda(|y|,|z|)$, and $w_\tau$ be given by (\ref{eqA.4}) with $\mu = \Lambda_0 - \lambda \rightarrow 0^{+}$. Define
$$
\widehat{\mathcal{N}}_{\tau} = \{w \in H: \frac{1}{\tau}\int_{\mathbb{R}^{3}}(|\nabla_y w|^2 + |y|^2w^2 - \Lambda_0w^2)dx + \int_{\mathbb{R}^{3}}(|\partial_z w|^{2} + w^{2})dx = \int_{\mathbb{R}^{3}}|w|^{p}dx\},
$$
$$
\widehat{\mathcal{N}}_{0} = \{w(z) \in H^1(\mathbb{R}): \int_{\mathbb{R}}(|w'(z)|^{2} + w(z)^{2})dx = \frac{2}{p}\pi^{1-\frac{p}{2}}\int_{\mathbb{R}}|w(z)|^{p}dx\},
$$
$$
\widehat{J}_0(w) = \frac{1}{2}\int_{\mathbb{R}}(|w'(z)|^{2} + w(z)^{2})dx - \frac{2}{p^2}\pi^{1-\frac{p}{2}}\int_{\mathbb{R}}|w(z)|^{p}dx
$$
and
$$
\widehat{h}(\tau) = \inf_{w \in \widehat{\mathcal{N}}_{\tau}}\widehat{J}_{\tau}(w).
$$

Similar to Lemma \ref{lem3.1}, we have

\begin{lemma}
\label{lemB.1}	$J_\lambda(u_\lambda) = h(\lambda) \Leftrightarrow \widehat{J}_{\tau}(w_\tau) = \widehat{h}(\tau)$ for any $\tau > 0$.
\end{lemma}

\begin{lemma}
\label{lemB.2} $\widehat{h}(\tau) \leq \widehat{h}(0), \forall \tau > 0$.
\end{lemma}

\textit{Proof.  } Note that $e_1(y)\widehat{w}(z) \in \mathcal{N}_{\tau}, \forall \tau > 0$. Hence, we have
$$
\widehat{h}(\tau) \leq \widehat{J}_{\tau}(e_1(y)\widehat{w}(z)) = \widehat{h}(0), \forall \tau > 0.
$$
\qed\vskip 5pt

\begin{lemma}
\label{lemB.3} $\int_{\mathbb{R}^{3}}(|\nabla_y w_\tau|^2 + |y|^2w_\tau^2 - \Lambda_0w_\tau^2)dx \lesssim \tau, \int_{\mathbb{R}^{3}}(|\partial_z w_\tau|^{2} + w_\tau^{2})dx \lesssim 1, \int_{\mathbb{R}^3}|w_\tau|^{p}dx \lesssim 1, \forall \tau > 0$.
\end{lemma}

\textit{Proof.  } Lemma \ref{lemB.2} yields that
$$
\frac{1}{2\tau}\int_{\mathbb{R}^{3}}(|\nabla_y w_\tau|^2 + |y|^2w_\tau^2 - \Lambda_0w_\tau^2)dx + \frac{1}{2}\int_{\mathbb{R}^{3}}(|\partial_z w_\tau|^{2} + w_\tau^{2})dx - \frac{1}{p}\int_{\mathbb{R}^3}|w_\tau|^{p}dx \leq \widehat{h}(0).
$$
Then, noticing
$$
\frac{1}{\tau}\int_{\mathbb{R}^{3}}(|\nabla_y w_\tau|^2 + |y|^2w_\tau^2 - \Lambda_0w_\tau^2)dx + \int_{\mathbb{R}^{3}}(|\partial_z w_\tau|^{2} + w_\tau^{2})dx = \int_{\mathbb{R}^{3}}|w_\tau|^{p}dx,
$$
we can deduce that
$$
(\frac{1}{2} - \frac{1}{p})\int_{\mathbb{R}^{3}}|w_\tau|^{p} \leq \widehat{h}(0),
$$
$$
(\frac{1}{2} - \frac{1}{p})(\frac{1}{\tau}\int_{\mathbb{R}^{3}}(|\nabla_y w_\tau|^2 + |y|^2w_\tau^2 - \Lambda_0w_\tau^2)dx + \int_{\mathbb{R}^{3}}(|\partial_z w_\tau|^{2} + w_\tau^{2})dx) \leq \widehat{h}(0).
$$
\qed\vskip 5pt

\begin{lemma}
\label{lemB.4} $w_\tau$ is uniformly bounded in $L^\infty(\mathbb{R}^3)$ and in $H^2_V(\mathbb{R}^3)$ with respect to $\tau \leq C$ for any $C > 0$, where
$$
H^2_V(\mathbb{R}^3) = \{w \in L^2(\mathbb{R}^3): (-\Delta + (x_1^2 + x_2^2))w \in L^2(\mathbb{R}^3)\}.
$$
\end{lemma}

\textit{Proof.  } From Lemma \ref{lemB.3} we derive that $w_\tau$ is uniformly bounded in $H$. Note that
$$
w_\tau = (\frac{1}{\tau}(-\Delta_y + |y|^2 - \Lambda_0) - \partial_{zz} + 1)^{-1}|w_\tau|^{p-2}w_\tau.
$$
Then, by iterating sufficiently many times, we conclude that $\|w_\tau\|_{L^\infty} \lesssim 1$. In particular, we know that $\|w_\tau\|_{L^{2(p-1)}} \lesssim 1$ in this iterating process. Hence,
\begin{eqnarray}
&&\|(-\Delta_y + |y|^2 -\partial_{zz})w_\tau\|_{L^2} \nonumber \\
&=& \|(-\Delta_y + |y|^2 -\partial_{zz})(\frac{1}{\tau}(-\Delta_y + |y|^2 - \Lambda_0) - \partial_{zz} + 1)^{-1}|w_\tau|^{p-2}w_\tau\|_{L^2} \nonumber \\
&\lesssim& \||w_\tau|^{p-2}w_\tau\|_{L^2} \nonumber \\
&=& \|w_\tau\|_{L^{2(p-1)}}^{p-1} \nonumber \\
&\lesssim& 1.
\end{eqnarray}
\qed\vskip 5pt

Let $\{e_i\}_{i = 1}^{\infty} \subset L^2_y(\mathbb{R}^2)$ be the collection of $L^2_y(\mathbb{R}^2)-$normalized eigenfunctions for $-\Delta_y + |y|^2$, i.e.
$$
(-\Delta_y + |y|^2)e_i = \lambda_ie_i,
$$
with eigenvalues $\Lambda_0 = \lambda_1 < \lambda_2 \leq \lambda_3 \leq \cdots$ in a non-decreasing order. Recall that $\{e_i\}_{i = 1}^{\infty} \subset L^2_y(\mathbb{R}^2)$ forms an orthonormal basis of $L^2_y(\mathbb{R}^2)$, and that the lowest eigenvalue $\Lambda_0 = 2$ is simple and the corresponding eigenfunction is given by $e_1 = \frac{1}{\sqrt{\pi}}e^{-\frac{|y|^2}{2}}$. From the spectral representation, $w \in L^2(\mathbb{R}^3)$ can be written as
$$
w(x) = \sum_{i = 1}^{\infty}\langle w(y,z),e_i(y)\rangle_{L^2_y(\mathbb{R}^2)}e_i(y), x = (y,z),
$$
where
$$
\langle w(y,z),e_i(y)\rangle_{L^2_y(\mathbb{R}^2)} = \int_{\mathbb{R}^2)}w(y,z),e_i(y)dy: \mathbb{R}_z \rightarrow \mathbb{R}.
$$
Denote $\langle w(y,z),e_1(y)\rangle_{L^2_y(\mathbb{R}^2)}$ by $\rho(z)$ and define
$$
(Pw)(x) = \rho(z)e_1(y): \mathbb{R}^3 \rightarrow \mathbb{R}.
$$
Let $Q = 1 - P$, i.e.
$$
(Qw)(x) = \sum_{i = 2}^{\infty}\langle w(y,z),e_i(y)\rangle_{L^2_y(\mathbb{R}^2)}e_i(y): \mathbb{R}^3 \rightarrow \mathbb{R}.
$$
From now on, we aim to show that $\rho_\tau(z) \rightarrow \widehat{w}(z)$ and $Qw_\tau \rightarrow 0$, where $\widehat{w}(z)$ is given by Theorem \ref{thmA.3} and $\rho_\tau(z) = \langle w_\tau(y,z),e_i(y)\rangle_{L^2_y(\mathbb{R}^2)}$. By doing that, the proof of Theorem \ref{thmA.3} would be complete.

Similar to Lemma \ref{lem3.2}, we have
\begin{lemma}
\label{lemB.5}	For any $\rho \in H^1_z(\mathbb{R}) \backslash \{0\}$, there exists a unique function $\widehat{t}_0: H^1_z(\mathbb{R}) \backslash \{0\} \rightarrow (0, +\infty)$ such that
\begin{equation}
t\rho \in \mathcal{\widehat{N}}_0, t > 0 \Leftrightarrow t = \widehat{t}_0(\rho),
\nonumber
\end{equation}
and $\widehat{J}_{0}(\widehat{t}_0(\rho)\rho) = \max_{t > 0}\widehat{J}_{0}(t\rho)$, $\widehat{t}_0 \in C(H^1_z(\mathbb{R}) \backslash \{0\}, (0, \infty))$.
\end{lemma}

\begin{lemma}
\label{lemB.6}  For any $C > 0$, there exists $\delta = \delta(C) > 0$ such that $\widehat{h}(\tau) \geq \delta$ for $\tau < C$.
\end{lemma}	

\textit{Proof.  }  For any $w \in H$ with $\|w\| = r$, we have
\begin{eqnarray}
\widehat{J}_\tau(w) \geq \frac{1}{2(C + \Lambda_0 + 1)}\|w\| - \frac{1}{p}\|w\|_{L^p(\mathbb{R}^3)}^p &\geq& \frac{1}{2(C + \Lambda_0 + 1)}r^2 - \widetilde{C}r^p,
\end{eqnarray}
where $\widetilde{C} > 0$ is a constant. Take $r$ small such that
$$
\delta = \frac{1}{2(C + \Lambda_0 + 1)}r^2 - \widetilde{C}r^p > 0.
$$
Thus the proof is completed.
\qed\vskip 5pt

\begin{theorem}
\label{thmB.7} $\rho_\tau(z) \rightarrow \widehat{w}$ in $H^1_z(\mathbb{R})$ as $\tau \rightarrow 0^+$.
\end{theorem}

\textit{Proof.  }  Step 1: $Qw_\tau \rightarrow 0$ in $L^q(\mathbb{R}^3), 2 \leq q < 2^\ast$ as $\tau \rightarrow 0^+$.

First, we show that $Qw_\tau \rightarrow 0$ in $L^2(\mathbb{R}^3)$. Indeed,
\begin{eqnarray}
\|Qw_\tau\|_{L^2(\mathbb{R}^3)}^2 &\lesssim& \int_{\mathbb{R}^{3}}(|\nabla_y (Qw_\tau)|^2 + |y|^2(Qw_\tau)^2)dx \nonumber \\
&\lesssim& \int_{\mathbb{R}^{3}}(|\nabla_y (Qw_\tau)|^2 + |y|^2(Qw_\tau)^2 - \Lambda_0(Qw_\tau)^2)dx \nonumber \\
&=& \int_{\mathbb{R}^{3}}(|\nabla_y w_\tau|^2 + |y|^2w_\tau^2 - \Lambda_0w_\tau^2)dx \nonumber \\
&\lesssim& \tau.
\end{eqnarray}
Then, noticing that $Qw_\tau$ is bounded in $H$, we can conclude.

Step 2: Let $w_n = w_{\tau_n}$, $\rho_n = \rho_{\tau_n}$ with $\tau_n \rightarrow 0^+$. Then up to a subsequence, $\rho_n(z) \rightarrow \widehat{w}$ in $H^1_z(\mathbb{R})$.

Lemma \ref{lemB.6} yields that $\|w_n\|_{L^p(\mathbb{R}^3)} \geq \epsilon > 0$. By \cite[Lemma 3.4]{BBJV}, from the boundedness of $w_n$ in $H$ and the fact that $w_n = w_n(|y|,|z|)$ is decreasing in $|y| > 0, |z| > 0$, we can conclude that up to a subsequence,
$$
w_n(y,z) \rightharpoonup w_0 \neq 0 \ in \ H, w_n(y,z) \rightarrow w_0 \neq 0 \ in \ L_{loc}^p(\mathbb{R}^3).
$$
By Step 1, we can verify that $w_0 = \rho_0(z)e_1(y)$ a.e. on $\mathbb{R}^3$, where $\rho_0(z) = \langle w_0(y,z),e_i(y)\rangle_{L^2_y(\mathbb{R}^2)}$.

For any $\varphi(z) \in C_0^\infty(\mathbb{R})$, $\varphi(z)e_1(y) \in H$. Thus, we have
\begin{eqnarray}
&& \int_{\mathbb{R}^3}(\rho_0'(z)e_1(y)\varphi'(z)e_1(y) + \rho_0(z)e_1(y)\varphi(z)e_1(y))dx \nonumber \\
&=& \lim_{n \rightarrow +\infty}\int_{\mathbb{R}^3}(\partial_zw_n(y,z)\varphi'(z)e_1(y) + w_n(y,z)\varphi(z)e_1(y))dx \nonumber \\
&=& \lim_{n \rightarrow +\infty}\int_{\mathbb{R}^3}|w_n(y,z)|^{p-2}w_n(y,z)\varphi(z)e_1(y)dx \nonumber \\
&=& \int_{\mathbb{R}^3}|\rho_0(z)|^{p-2}\rho_0(z)\varphi(z)e_1(y)^{p-1}dx.
\end{eqnarray}
Since $C_0^\infty(\mathbb{R})$ is dense in $H_z^1(\mathbb{R})$, we deduce that
$$
\widehat{J}_0'(\rho_0(z)) = 0  \ in \ H_z^{-1}(\mathbb{R}),
$$
yielding that $\rho_0(z) \in \mathcal{\widehat{N}}_0$. Then, from
\begin{eqnarray}
\widehat{h}(0) &\leq& \widehat{J}_0(\rho_0) \nonumber \\
&=& \frac{p - 2}{p^2}\pi^{1-\frac{p}{2}}\int_{\mathbb{R}}|\rho_0|^{p}dz \nonumber \\
&\leq& \lim_{n \rightarrow +\infty}\frac{p - 2}{p^2}\pi^{1-\frac{p}{2}}\int_{\mathbb{R}}|\rho_n(z)|^{p}dz \nonumber \\
&=& \lim_{n \rightarrow +\infty}\frac{p - 2}{2p}\int_{\mathbb{R}^3}|\rho_n(z)e_1(y)|^{p}dx \nonumber \\
&=& \lim_{n \rightarrow +\infty}\frac{p - 2}{2p}\int_{\mathbb{R}^3}|w_n|^{p}dx \nonumber \\
&=& \lim_{n \rightarrow +\infty}\widehat{h}(\tau_n) \nonumber \\
&\leq& \widehat{h}(0),
\end{eqnarray}
we obtain that $\widehat{J}_0(\rho_0) = \widehat{h}(0)$, yielding that $\rho_0 = \widehat{w}$, and that
$$
\int_{\mathbb{R}}|\rho_n(z)|^{p}dz \rightarrow \int_{\mathbb{R}}|\rho_0|^{p}dz,
$$
implying that
$$
\int_{\mathbb{R}}(|\rho_n'(z)|^2 + \rho_n(z)^2)dz \rightarrow \int_{\mathbb{R}}(|\rho_0'|^2 + \rho_0^2)dz.
$$
Noticing the weak convergence of $\rho_n(z)$, we know that $\rho_n(z) \rightarrow \widehat{w}$ in $H^1_z(\mathbb{R})$. This completes the proof.
\qed\vskip 5pt

Define
$$
\widehat{L}_\tau = \frac{1}{\tau}(-\Delta_y + |y|^2 - \Lambda_0) - \partial_{zz} + 1 - (p-1)|w_\tau|^{p-2}, \tau > 0,
$$
which is a self-adjoint operator from $L^2(\mathbb{R}^3)$ to $L^2(\mathbb{R}^3)$, with form domain $H$. We say that $w_\tau$ is non-degenerate in $H$ if $\ker \widehat{L}_\tau = \{\partial_z w_\tau\}$.

\textbf{Proof of Theorem \ref{thmA.3}.  } By Theorem \ref{thmB.7}, it suffices to show that
$$
Qw_\tau \rightarrow 0 \ in \ H.
$$
By Lemma \ref{lemB.5} and Theorem \ref{thmB.7}, $\widehat{t}_0(\rho_\tau(z)) \rightarrow 1$. Then
\begin{eqnarray}
\widehat{J}_0(\widehat{w}) &\leq& \widehat{J}_0(\widehat{t}_0(\rho_\tau(z))\rho_\tau(z)) \nonumber \\
&=& \widehat{J}_0(\rho_\tau(z)) + o(1) \nonumber \\
&=& \widehat{J}_\tau(\rho_\tau(z)e_1(y)) + o(1) \nonumber \\
&=& \widehat{J}_\tau(w_\tau) - \frac{1}{2\tau}\int_{\mathbb{R}^{3}}(|\nabla_y (Qw_\tau)|^2 + |y|^2(Qw_\tau)^2 - \Lambda_0(Qw_\tau)^2)dx \nonumber \\
&& - \frac{1}{2}\int_{\mathbb{R}^{3}}(|\partial_z (Qw_\tau)|^{2} + (Qw_\tau)^{2})dx + o(1).
\end{eqnarray}
Using Lemma \ref{lemB.2}, $\widehat{J}_0(\widehat{w}) \geq \widehat{J}_\tau(w_\tau)$, we obtain that:
$$
\frac{1}{\tau}\int_{\mathbb{R}^{3}}(|\nabla_y (Qw_\tau)|^2 + |y|^2(Qw_\tau)^2 - \Lambda_0(Qw_\tau)^2)dx \rightarrow 0 \ as \ \tau \rightarrow 0^+.
$$
$$
\int_{\mathbb{R}^{3}}(|\partial_z (Qw_\tau)|^{2} + (Qw_\tau)^{2})dx \rightarrow 0 \ as \ \tau \rightarrow 0^+.
$$
Thus we show the convergence of $w_\tau$ in $H$.

Then, along the lines of \cite[Section 5]{HJ}, we can show the non-degeneracy and uniqueness of $w_\tau$ in $H$ for sufficiently small $\tau$, which is equivalent to the non-degeneracy and uniqueness of $u_\lambda$ in $H$ with $\Lambda_2 < \lambda < \Lambda_0$ for some $\Lambda_2 < \Lambda_0$.

Finally, noticing that:

\begin{eqnarray}
&& (\Lambda_0 - \lambda)^{\frac{1}{2} - \frac{2}{p-2}}\int_{\mathbb{R}^{3}}(|\nabla_y u_\lambda|^2 + u_\lambda^2)dx + (\Lambda_0 - \lambda)^{-\frac{1}{2} - \frac{2}{p-2}}\int_{\mathbb{R}^{3}}|\partial_z u_\lambda|^{2}dx \nonumber \\
&=& \int_{\mathbb{R}^{3}}(|\nabla_y w_\tau|^2 + w_\tau^2 + |\partial_z w_\tau|^{2})dx \nonumber \\
&\rightarrow& \|e_1(y)\widehat{w}(z)\|_{H^1(\mathbb{R}^{3})}^2,
\end{eqnarray}
and
\begin{eqnarray}
&& (\Lambda_0 - \lambda)^{\frac{3}{2} - \frac{2}{p-2}}\int_{\mathbb{R}^{3}}(|\nabla_y u_\lambda|^2 + u_\lambda^2)dx \nonumber \\
&=& \int_{\mathbb{R}^{3}}|y|^2w_\tau^2dx \nonumber \\
&\rightarrow& \int_{\mathbb{R}^2}|y|^2e_1(y)^2dy\int_{\mathbb{R}}\widehat{w}(z)^2dz,
\end{eqnarray}
we can obtain the asymptotic behaviors of $\|u_\lambda\|_{H^1(\mathbb{R}^3)}$ and $\|u_\lambda\|$ as $\lambda \rightarrow \Lambda_0$. This completes the proof.
\qed\vskip 5pt

\section{The existence, non-existence, multiplicity of normalized solutions, and orbital stability/instability }

\subsection{Asymptotical behavior and monotonicity of $\int_{\mathbb{R}^{3}}u_\lambda^2dx$ when $\lambda \rightarrow -\infty$}

Similar arguments to the proof of \cite[Corollary 2.5]{Song} and \cite[Lemma 4.23]{Song3}, yield that $\{(\lambda,u_{\lambda}): \lambda < \Lambda_{1}\}$ is a $C^{1}$ curve in $\mathbb{R} \times H$. Let
$$
\chi_\lambda = \partial_\lambda u_\lambda,
$$
$$
\widetilde{\chi}_\mu = |\lambda|^{-\frac{3-p}{p-2}}\chi_\lambda(\frac{x}{\sqrt{|\lambda|}})
$$
where $\mu = 1/\lambda^2$.

\begin{lemma}
\label{lem4.1}	Recall that $\widetilde{L}_{\mu,s}$ is given by Section 3.

$(i)$ For $\mu <  1/\Lambda_1^2$, there exists $C > 0$ such that
$$
\limsup_{\mu \rightarrow 0^+}\|\widetilde{L}_{\mu,s}^{-1}\|_{L^2(\mathbb{R}^3) \rightarrow H^2(\mathbb{R}^3)} \leq C,
$$
	
$(ii)$ For any $v \in L_s^2(\mathbb{R}^3)$, $\widetilde{L}_{\mu,s}^{-1}v \rightarrow \widetilde{L}_{0,s}^{-1}v$ in $H^2(\mathbb{R}^3)$.
\end{lemma}

\textit{Proof.  }  $(i)$ Note that $\widetilde{L}_\mu$ is $\widetilde{J}''_{\mu}(\widetilde{v})$. By the proof of Theorem \ref{thm1.3}, we may assume that there exists $\epsilon > 0$ such that
$$
\sigma(\widetilde{L}_{\mu,s}) \cap (-\epsilon,\epsilon) = \emptyset, \forall \mu \in [0,1/\Lambda_1^2).
$$
Hence,
$$
\|\widetilde{L}_{\mu,s}^{-1}v\|_{H^2(\mathbb{R}^3)} \leq \|\widetilde{L}_{\mu,s}^{-1}v\|_{H_V^2(\mathbb{R}^3)} \leq \frac{1}{\epsilon}\|v\|_{L^2(\mathbb{R}^3)}.
$$
Where
$$
H_V^2(\mathbb{R}^3) = \{v \in L^2(\mathbb{R}^N): (-\Delta v + (x_1^2 + x_2^2)v) \in L^2(\mathbb{R}^N)\}
$$
is the domain of $\widetilde{L}_\mu$.

$(ii)$ By the proof of $(i)$, $\widetilde{L}_{\mu,s}^{-1}v$ is bounded in $H^2_V$ for any $v \in L_s^2(\mathbb{R}^3)$. Therefore,
$$
\widetilde{L}_{0,s}\widetilde{L}_{\mu,s}^{-1}v = v - \mu(x_1^2+x_2^2)\widetilde{L}_{\mu,s}^{-1}v \rightarrow v \in L^2(\mathbb{R}^3),
$$
implying that
$$
\|\widetilde{L}_{\mu,s}^{-1}v - \widetilde{L}_{0,s}^{-1}v\|_{H^2(\mathbb{R}^3)} = \|\widetilde{L}_{0,s}^{-1}(\widetilde{L}_{0,s}\widetilde{L}_{\mu,s}^{-1}v - v)\|_{H^2(\mathbb{R}^3)} \leq C\|\widetilde{L}_{0,s}\widetilde{L}_{\mu,s}^{-1}v - v\|_{L^2(\mathbb{R}^3)} \rightarrow 0.
$$
\qed\vskip 5pt

\begin{lemma}

\label{lem4.2}  Let $10/3 < p < 6$. Then

$(i)$ $\int_{\mathbb{R}^{3}}u_\lambda^2dx \rightarrow 0$ as $\lambda \rightarrow -\infty$;

$(ii)$ $\widetilde{\chi}_\mu \rightarrow \frac{1}{2-p}\widetilde{v} - \frac{1}{2}x \cdotp \nabla \widetilde{v}$ in $H^2(\mathbb{R}^3)$ as $\mu \rightarrow 0^+$.

$(iii)$ for some $\widetilde{\Lambda}_1 \leq \Lambda_1$, $\frac{d}{d\lambda} \int_{\mathbb{R}^{3}}u_\lambda^2dx > 0$ when $\lambda < \widetilde{\Lambda}_1$.
\end{lemma}

\textit{Proof.  }  $(i)$ From Theorem \ref{thm1.3}, we derive that
$$
|\lambda|^{\frac{3}{2}-\frac{2}{p-2}}\int_{\mathbb{R}^{3}}u_\lambda^2dx \rightarrow \int_{\mathbb{R}^{3}}\widetilde{v}^2dx,
$$
showing that $\int_{\mathbb{R}^{3}}u_\lambda^2dx \rightarrow 0$, since $3/2 - 2/(p-2) > 0$.

$(ii)$ Direct computations show that
$$
\widetilde{L}_{\mu,s}\widetilde{\chi}_\mu = v_\mu.
$$
Together with Lemma \ref{lem4.1} and the convergence of $v_\mu$, we have
$$
\widetilde{\chi}_\mu = \widetilde{L}_{\mu,s}^{-1}v_\mu = \widetilde{L}_{\mu,s}^{-1}(v_\mu - \widetilde{v}) + \widetilde{L}_{\mu,s}^{-1}\widetilde{v} \rightarrow \widetilde{L}_{0,s}^{-1}\widetilde{v} \ in \ H^2(\mathbb{R}^3).
$$
Noticing that
$$
\widetilde{L}_{0,s}(\frac{1}{2-p}\widetilde{v} - \frac{1}{2}x \cdotp \nabla \widetilde{v}) = \widetilde{v},
$$
we can conclude.

$(iii)$ \begin{eqnarray}
&& \lim_{\lambda \rightarrow -\infty}|\lambda|^{\frac{p-4}{p-2} + \frac{3}{2}}\frac{d}{d\lambda}\int_{\mathbb{R}^{3}}|u_{\lambda}|^{2}dx \nonumber \\
&=& \lim_{\lambda \rightarrow -\infty}2|\lambda|^{\frac{p-4}{p-2} + \frac{3}{2}}\int_{\mathbb{R}^{N}}u_{\lambda}\chi_{\lambda}dx \nonumber \\
&=& \lim_{\mu \rightarrow 0^+}2\int_{\mathbb{R}^{3}}v_{\mu}\widetilde{\chi}_{\mu}dx \nonumber \\
&=& 2\int_{\mathbb{R}^{3}}\widetilde{v}(\frac{1}{2-p}\widetilde{v} - \frac{1}{2}x \cdotp \nabla \widetilde{v})dx \nonumber \\
&=& 2(\frac{1}{2-p}+\frac{3}{4})\int_{\mathbb{R}^{3}}|\widetilde{v}|^{2}dx > 0.
\end{eqnarray}
Therefore, there exists $\widetilde{\Lambda}_{1} \leq \Lambda_{1}$ such that for any $\lambda < \widetilde{\Lambda}_{1}$,
$$
\frac{d}{d\lambda}\int_{\mathbb{R}^{3}}|u_{\lambda}|^{2}dx > 0.
$$
\qed\vskip 5pt

\subsection{Asymptotical behavior and monotonicity of $\int_{\mathbb{R}^{3}}u_\lambda^2dx$ when $\lambda \rightarrow \Lambda_{0}$}

Similar arguments to the proof of \cite[Corollary 2.5]{Song} and \cite[Lemma 4.23]{Song3}, yield that $\{(\lambda,u_{\lambda}): \Lambda_{2} < \lambda < \Lambda_{0}\}$ is a $C^{1}$ curve in $\mathbb{R} \times H$. Let
$$
\chi_\lambda = \partial_\lambda u_\lambda,
$$
$$
\widehat{\chi}_\tau = (\Lambda_{0} - \lambda)^{-\frac{3-p}{p-2}}\chi_\lambda(y,\frac{z}{\sqrt{\Lambda_{0} - \lambda}})
$$
where $\tau = \Lambda_{0} - \lambda$.

\begin{lemma}
\label{lem4.3}	Recall that $\widehat{L}_\tau$ is given by Section 4 and define
$$
\widehat{L}_0 = - \frac{d^2}{dz^2} + 1 - \frac{2(p-1)}{p}\pi^{1-\frac{p}{2}}|\widehat{w}(z)|^{p-2},
$$
which is a self-adjoint operator from $L^2(\mathbb{R})$ to $L^2(\mathbb{R})$, with form domain $H^1(\mathbb{R})$. Then there exists $\tau_0 > 0$ and $\epsilon > 0$ such that
$$
\|\widehat{L}_\tau\phi\|_{L^2(\mathbb{R}^3)} \geq \epsilon\|\phi\|
$$
for any $\phi = \phi(|y|,|z|) \in H^2_V$, $\tau \in (0,\tau_0)$.
\end{lemma}

\textit{Proof.  } Suppose on the contrary that there exist $\tau_n \rightarrow 0^+$ and $\phi_n = \phi_n(|y|,|z|) \in H^2_V$ such that $\|\phi_n\| = 1$ and
$$
\|\widehat{L}_{\tau_n}\phi_n\|_{L^2(\mathbb{R}^3)} \rightarrow 0.
$$
Passing to a subsequence if necessary, we may assume that $\phi_n \rightharpoonup \phi_\infty$ in $H$ and $\phi_n \rightarrow \phi_\infty$ in $L_{loc}^p(\mathbb{R}^3)$. Similar to the proof of Theorem \ref{thmB.7}, $\phi_\infty = \rho_\infty(z)e_1(y)$ a.e. on $\mathbb{R}^3$, where $\rho_\infty(z) = \langle \phi_\infty(y,z),e_i(y)\rangle_{L^2_y(\mathbb{R}^2)}$ and for any $\varphi(z) \in C_0^\infty(\mathbb{R})$, we have
\begin{eqnarray}
&& \int_{\mathbb{R}^3}(\rho_\infty'(z)e_1(y)\varphi'(z)e_1(y) + \rho_\infty(z)e_1(y)\varphi(z)e_1(y))dx \nonumber \\
&=& \lim_{n \rightarrow +\infty}\int_{\mathbb{R}^3}(\partial_z\phi_n(y,z)\varphi'(z)e_1(y) + \phi_n(y,z)\varphi(z)e_1(y))dx \nonumber \\
&=& \lim_{n \rightarrow +\infty}(\langle \widehat{L}_{\tau_n}\phi_n,\varphi(z)e_1(y) \rangle_{L^2(\mathbb{R}^3)} + (p-1)\int_{\mathbb{R}^3}|w_{\tau_n}|^{p-2}\phi_n\varphi(z)e_1(y)dx) \nonumber \\
&=& (p-1)\int_{\mathbb{R}^3}|\widehat{w}(z)|^{p-2}\rho_\infty(z)\varphi(z)e_1(y)^{p}dx.
\end{eqnarray}
Then the non-degeneracy of $\widehat{w}(z)$ yields that $\rho_\infty(z) = 0$, and hence, $\phi_\infty = 0$. Note that
$$
1 = \|\phi_n\| \lesssim \langle \widehat{L}_{\tau_n}\phi_n,\phi_n \rangle_{L^2(\mathbb{R}^3)} + (p-1)\int_{\mathbb{R}^3}|w_{\tau_n}|^{p-2}\phi_n^2dx.
$$
However, form the convergence of $w_{\tau_n}$ and $\phi_n$, we obtain that
$$
\int_{\mathbb{R}^3}|w_{\tau_n}|^{p-2}\phi_n^2dx \rightarrow 0,
$$
implying that
$$
\langle \widehat{L}_{\tau_n}\phi_n,\phi_n \rangle_{L^2(\mathbb{R}^3)} + (p-1)\int_{\mathbb{R}^3}|w_{\tau_n}|^{p-2}\phi_n^2dx \rightarrow 0,
$$
which is absurd and we complete the proof.
\qed\vskip 5pt

\begin{lemma}

\label{lem4.4}  Let $2 < p < 6$. Then

$(i)$ $\int_{\mathbb{R}^{3}}u_\lambda^2dx \rightarrow 0$ as $\lambda \rightarrow \Lambda_0$;

$(ii)$ $\widehat{\chi}_\tau \rightarrow (\frac{1}{2-p}\widehat{w}(z) - \frac{1}{2}z\widehat{w}'(z))e_1(y)$ in $H$ as $\tau \rightarrow 0^+$.

$(iii)$ for some $\widetilde{\Lambda}_2 \in [\Lambda_2,\Lambda_0)$, $\frac{d}{d\lambda} \int_{\mathbb{R}^{3}}u_\lambda^2dx < 0$ when $\lambda \in (\widetilde{\Lambda}_2,\Lambda_0)$.
\end{lemma}

\textit{Proof.  }  $(i)$ From Theorem \ref{thmA.3}, we derive that
$$
(\Lambda_0 - \lambda)^{\frac{1}{2}-\frac{2}{p-2}}\int_{\mathbb{R}^{3}}u_\lambda^2dx \rightarrow \int_{\mathbb{R}}\widehat{w}^2dz,
$$
showing that $\int_{\mathbb{R}^{3}}u_\lambda^2dx \rightarrow 0$, since $1/2 - 2/(p-2) < 0$.

$(ii)$ Direct computations show that
$$
\widehat{L}_{\tau}\widehat{\chi}_\tau = w_\tau.
$$
By Lemma \ref{lem4.3},
$$
\|\widehat{\chi}_\tau\| \leq \frac{1}{\epsilon}\|\widehat{L}_{\tau}\widehat{\chi}_\tau\|_{L^2(\mathbb{R}^3)} = \frac{1}{\epsilon}\|w_\tau\|_{L^2(\mathbb{R}^3)} \lesssim 1.
$$
For any $\tau_n \rightarrow 0$, up to a subsequence, we may assume that $\widehat{\chi}_{\tau_n} \rightharpoonup \widehat{\chi}_0$ in $H$ and $\widehat{\chi}_{\tau_n} \rightarrow \widehat{\chi}_0$ in $L_{loc}^q(\mathbb{R}^3), 2 \leq q < 2^\ast$. Standard arguments yields that $\widehat{\chi}_0 = \varrho_0(z)e_1(y)$, where $\varrho_0(z) = \langle \widehat{\chi}_0(y,z),e_i(y)\rangle_{L^2_y(\mathbb{R}^2)}$, and $\widehat{L}_0\varrho_0(z) = \widehat{w}(z)$. Noticing that $\widehat{L}_0$ is invertible and
$$
\widehat{L}_0(\frac{1}{2-p}\widehat{w} - \frac{1}{2}z\widehat{w}') = \widehat{w},
$$
we obtain that $\varrho_0 = \frac{1}{2-p}\widehat{w} - \frac{1}{2}z\widehat{w}'$.

Similar to Step 1 in the proof of Theorem \ref{thmB.7}, we have $Q\widehat{\chi}_{\tau_n} \rightarrow 0$ in $L^q(\mathbb{R}^3), 2 \leq q < 2^\ast$. Then,
\begin{eqnarray}
&& \frac{1}{\tau_n}\int_{\mathbb{R}^{3}}(|\nabla_y (Q\widehat{\chi}_{\tau_n})|^2 + |y|^2(Q\widehat{\chi}_{\tau_n})^2 - \Lambda_0(Q\widehat{\chi}_{\tau_n})^2)dx + \int_{\mathbb{R}^{3}}(|\partial_z (Q\widehat{\chi}_{\tau_n})|^{2} + (Q\widehat{\chi}_{\tau_n})^{2})dx \nonumber \\
&=& \int_{\mathbb{R}^{3}}(Q\widehat{\chi}_{\tau_n})^2dx + (p-1) \int_{\mathbb{R}^{3}}|w_{\tau_n}|^{p-2}\widehat{\chi}_{\tau_n}Q\widehat{\chi}_{\tau_n}dx \nonumber \\
&\rightarrow& 0,
\end{eqnarray}
implying that
\begin{eqnarray}
&& \frac{1}{\tau_n}\int_{\mathbb{R}^{3}}(|\nabla_y (Q\widehat{\chi}_{\tau_n})|^2 + |y|^2(Q\widehat{\chi}_{\tau_n})^2)dx \nonumber \\
&\lesssim& \frac{1}{\tau_n}\int_{\mathbb{R}^{3}}(|\nabla_y (Q\widehat{\chi}_{\tau_n})|^2 + |y|^2(Q\widehat{\chi}_{\tau_n})^2 - \Lambda_0(Q\widehat{\chi}_{\tau_n})^2)dx \nonumber \\
&\rightarrow& 0,
\end{eqnarray}
and
$$
\int_{\mathbb{R}^{3}}(|\partial_z (Q\widehat{\chi}_{\tau_n})|^{2} + (Q\widehat{\chi}_{\tau_n})^{2})dx \rightarrow 0.
$$
Thus $Q\widehat{\chi}_{\tau_n} \rightarrow 0$ in $H$. Moreover,
\begin{eqnarray}
\|\varrho_{\tau_n}\|_{H^1(\mathbb{R})}^2 &=& \int_{\mathbb{R}^{3}}(|\varrho'_{\tau_n}(z)|^{2}e_1(y)^2 + \varrho_{\tau_n}^{2}e_1(y)^2)dx \nonumber \\
&=& \langle \widehat{L}_{\tau_n}\widehat{\chi}_{\tau_n},P\widehat{\chi}_{\tau_n} \rangle_{L^2(\mathbb{R}^3)} + (p-1) \int_{\mathbb{R}^{3}}|w_{\tau_n}|^{p-2}(P\widehat{\chi}_{\tau_n})\widehat{\chi}_{\tau_n}dx \nonumber \\
&=& \int_{\mathbb{R}^{3}}w_{\tau_n}(P\widehat{\chi}_{\tau_n})dx + (p-1) \int_{\mathbb{R}^{3}}|w_{\tau_n}|^{p-2}(P\widehat{\chi}_{\tau_n})\widehat{\chi}_{\tau_n}dx \nonumber \\
&\rightarrow& \int_{\mathbb{R}}\widehat{w}(z)\varrho_0(z)dz + \frac{2(p-1)}{p}\pi^{1-\frac{p}{2}}\int_{\mathbb{R}}|\widehat{w}(z)|^{p-2}\varrho_0(z)^2dz \nonumber \\
&=& \langle \widehat{L}_{0}\varrho_{0},\varrho_0 \rangle_{L^2(\mathbb{R})} + \frac{2(p-1)}{p}\pi^{1-\frac{p}{2}}\int_{\mathbb{R}}|\widehat{w}(z)|^{p-2}\varrho_0(z)^2dz \nonumber \\
&=& \|\varrho_{0}\|_{H^1(\mathbb{R})}^2,
\end{eqnarray}
implying that $\varrho_{\tau_n} \rightarrow \varrho_{0}$ in $H^1(\mathbb{R})$. Hence, $\widehat{\chi}_{\tau_n} \rightarrow \varrho_{0}(z)e_1(y)$ in $H$ and we can conclude.

$(iii)$ \begin{eqnarray}
&& \lim_{\lambda \rightarrow \Lambda_0}(\Lambda_0 - \lambda)^{\frac{p-4}{p-2} + \frac{1}{2}}\frac{d}{d\lambda}\int_{\mathbb{R}^{3}}|u_{\lambda}|^{2}dx \nonumber \\
&=& \lim_{\lambda \rightarrow \Lambda_0}2(\Lambda_0 - \lambda)^{\frac{p-4}{p-2} + \frac{1}{2}}\int_{\mathbb{R}^{N}}u_{\lambda}\chi_{\lambda}dx \nonumber \\
&=& \lim_{\tau \rightarrow 0^+}2\int_{\mathbb{R}^{3}}w_{\tau}\widehat{\chi}_{\tau}dx \nonumber \\
&=& 2\int_{\mathbb{R}}\widehat{w}(\frac{1}{2-p}\widehat{w} - \frac{1}{2}z \widehat{w}')dz \nonumber \\
&=& 2(\frac{1}{2-p}+\frac{1}{4})\int_{\mathbb{R}}|\widehat{w}|^{2}dz < 0.
\end{eqnarray}
Therefore, there exists $\widetilde{\Lambda}_2 \in [\Lambda_2,\Lambda_0)$ such that for any $\lambda \in (\widetilde{\Lambda}_2,\Lambda_0)$,
$$
\frac{d}{d\lambda}\int_{\mathbb{R}^{3}}|u_{\lambda}|^{2}dx < 0.
$$
\qed\vskip 5pt

\subsection{Proof of Theorem \ref{thm1.4}}

Recall that standing waves are called orbitally stable if for each $\epsilon > 0$ there
exists $\delta > 0$ such that, whenever $\Psi_0 \in H^1_V(\mathbb{R}^N , \mathbb{C})$ is such that $\|\Psi_0 - u_\lambda\|_{H^1_V(\mathbb{R}^N , \mathbb{C})} < \delta$ and $\Psi(t,x)$is the solution of (\ref{eq1.1}) with $\Psi(0,x) = \Psi_0(x)$ in some interval $[0, t_0)$, then $\Psi(t,x)$ can be continued to a solution in $0 \leq t < +\infty$ and
$$
\sup_{t > 0}\inf_{w \in \mathbb{R}, z \in \mathbb{R}}\|\Psi(t,x) - e^{-i\lambda w}u_\lambda(x_1,x_2,x_3-z)\|_{H^1_V(\mathbb{R}^N , \mathbb{C})} < \epsilon;
$$
otherwise, they are called unstable. To study the orbital stability, we lean on the following result, which expresses in our context the abstract theory developed in \cite{GSS2}.

Global well-posedness of the Cauchy problem was established in \cite{ACS}.

\begin{proposition}
\label{prop4.2} If $\frac{d}{d\lambda}\int_{\mathbb{R}^3}u_{\lambda}^{2} > 0$ (respectively $< 0$), the standing wave $e^{-i\lambda t}u_{\lambda}(x)$ is orbitally unstable (respectively stable) in
$$
H_V^1(\mathbb{R}^3,\mathbb{C}) = \{v \in H^1(\mathbb{R}^3,\mathbb{C}): \int_{\mathbb{R}^3}(x_1^2+x_2^2)|v|^2 dx < +\infty\}.
$$
\end{proposition}

\textbf{Proof of Theorem \ref{thm1.4}.  } Take
$$
\widehat{c} = \min\{\sqrt{\int_{\mathbb{R}^3}u_\lambda^2dx}: \lambda = \widetilde{\Lambda}_1 \ or \ \lambda = \widetilde{\Lambda}_2\},
$$
where $\widetilde{\Lambda}_1$ and $\widetilde{\Lambda}_2$ are given by Lemmas \ref{lem4.2} and \ref{lem4.4} respectively. Similar arguments to the proof of \cite[Corollary2.5]{Song} and \cite[Lemma 4.23]{Song3}, yield that
$$
\lambda \mapsto u_\lambda, \lambda < \widetilde{\Lambda}_1 \ or \ \widetilde{\Lambda}_2 < \lambda < \Lambda_0
$$
is $C^{1}$in $H$, and hence in $L^2(\mathbb{R}^3)$. Thus (\ref{eq1.2}) admits a ground state $u_\lambda$ with $\lambda < \widetilde{\Lambda}_1$ and a ground state $u_{\widetilde{\lambda}}$ with $\widetilde{\Lambda}_2 < \widetilde{\lambda} < \Lambda_0$, where
$$
\sqrt{\int_{\mathbb{R}^3}u_\lambda^2dx} = \sqrt{\int_{\mathbb{R}^3}u_{\widetilde{\lambda}}^2dx} = c, 0 < c < \widehat{c}.
$$

Note that the local well posedness in $H_V^1(\mathbb{R}^3,\mathbb{C})$ is established in \cite{ACS}. Then Lemma \ref{lem4.2} $(iii)$, Lemma \ref{lem4.4} $(iii)$, and Proposition \ref{prop4.2} yield the orbital stability of $e^{-i\widetilde{\lambda} t}u_{\widetilde{\lambda}}(x)$ in $H_V^1(\mathbb{R}^3,\mathbb{C})$, and the orbital instability of $e^{-i\lambda t}u_{\lambda}(x)$ in $H_{V,s}^1(\mathbb{R}^3,\mathbb{C})$, where
$$
H_{V,s}^1(\mathbb{R}^3,\mathbb{C}) = \{v \in H_V^1(\mathbb{R}^3,\mathbb{C}): v = v(|y|,|z|)\}.
$$
\qed\vskip 5pt

\subsection{Proof of Theorem \ref{thmA.8}}

We show some lemmas first. These results can be proved along the lines of Lemmas \ref{lem3.2}, \ref{lem3.3} and we omit the details here.

\begin{lemma}
For any $(\lambda, u) \in (-\infty,\Lambda_{0}) \times H \backslash \{0\}$, there exists a unique function $t: (-\infty,\Lambda_{0}) \times H \backslash \{0\} \rightarrow (0, +\infty)$ such that
\begin{equation}	
(\mu, t v) \in \mathcal{N}, t > 0 \Leftrightarrow t = t(\lambda, u),
\nonumber
\end{equation}
and $J_{\lambda}(t(\lambda, u)u) = \max_{t > 0}J_{\lambda}(tu)$, $t \in C((-\infty,\Lambda_{0}) \times H \backslash \{0\}, (0, \infty))$, where $\mathcal{N} := \cup_{\lambda < \Lambda_{0}}\mathcal{N}_{\lambda}$.
\end{lemma}

\begin{lemma}
\label{lem4.7}	For any $\lambda_0 < \Lambda_{0}$, $\limsup_{\lambda \rightarrow \lambda_0}h(\lambda) \leq h(\lambda_0)$.
\end{lemma}

\textbf{Proof of Theorem \ref{thmA.8}.  }  We will show that
\begin{equation} \label{eq5.7}
\sup_{\lambda \in [\widetilde{\Lambda}_1,\widetilde{\Lambda}_2]}\|u_\lambda\|_{L^2} \leq C,
\end{equation}
for some $C > 0$ independent of $\lambda$ where $\widetilde{\Lambda}_1, \widetilde{\Lambda}_2$ are given by Lemmas \ref{lem4.2}, \ref{lem4.4} respectively. This together with results proved in Lemmas \ref{lem4.2}, \ref{lem4.4}, will complete the proof.

By Lemma \ref{lem4.7}, there exists $\widetilde{C} > 0$ independent of $\lambda$ such that $\sup_{\lambda \in [\widetilde{\Lambda}_1,\widetilde{\Lambda}_2]} h(\lambda) \leq \widetilde{C}$. Since $J_\lambda(u_\lambda) = h(\lambda)$, we have
\begin{equation} \label{eq5.8}
\frac{1}{2}\int_{\mathbb{R}^{3}}(|\nabla u_\lambda|^{2} + (x_{1}^{2} + x_{2}^{2} - \lambda)u_\lambda^{2})dx - \frac{1}{p}\int_{\mathbb{R}^{3}}|u_\lambda|^{p}dx \leq \widetilde{C}.
\end{equation}
Moreover, $u_\lambda \in \mathcal{N}_{\lambda}$. Hence,
\begin{equation} \label{eq5.9}
\int_{\mathbb{R}^{3}}(|\nabla u_\lambda|^{2} + (x_{1}^{2} + x_{2}^{2} - \lambda)u_\lambda^{2})dx = \int_{\mathbb{R}^{3}}|u_\lambda|^{p}dx.
\end{equation}
From (\ref{eq5.8}) and (\ref{eq5.9}), we derive that
\begin{equation}
(\frac{1}{2} - \frac{1}{p})\int_{\mathbb{R}^{3}}(|\nabla u_\lambda|^{2} + (x_{1}^{2} + x_{2}^{2} - \lambda)u_\lambda^{2})dx \leq \widetilde{C}.
\end{equation}
Thus, we have
\begin{equation}
(\Lambda_{0} - \widetilde{\Lambda}_2)\sup_{\lambda \in [\widetilde{\Lambda}_1,\widetilde{\Lambda}_2]}\int_{\mathbb{R}^{3}}u_\lambda^{2}dx \leq  \sup_{\lambda \in [\widetilde{\Lambda}_1,\widetilde{\Lambda}_2]}\int_{\mathbb{R}^{3}}(|\nabla u_\lambda|^{2} + (x_{1}^{2} + x_{2}^{2} - \lambda)u_\lambda^{2})dx \leq \frac{p-2}{2p}\widetilde{C}.
\end{equation}
Take
$$
C = \sqrt{\frac{p-2}{2p(\Lambda_{0} - \widetilde{\Lambda}_2)}\widetilde{C}},
$$
and we show (\ref{eq5.7}). Thus we arrive at our conclusion.
\qed\vskip 5pt

\textbf{Acknowledgement:} The authors would like to thank Louis Jeanjean for providing valuable feedback and important references, and thank the Photonics group at Caltech for useful discussions and insight. The first author thanks C. Li and S.J. Li for important discussions.

\end{document}